\documentclass[aos,MSNbibl,seceqn,nameyear,dvips]{arximspdf}
\usepackage{mathrsfs}
\usepackage{dcolumn}

\doi{10.1214/11-AOS922}
\volume{39}
\issue{6}
\pubyear{2011}
\firstpage{2820}
\lastpage{2851}

\makeatletter

\newcolumntype{d}[1]{D{.}{.}{#1}}

\newcommand{\fracb}[2]{(#1)(#2)^{-1}}

\newtheorem{tm}{Theorem}
\newtheorem{la}{Lemma}

\newtheorem{pn}{Proposition}

\makeatother

\begin{document}
\begin{frontmatter}

\title{On the approximate maximum likelihood estimation for diffusion processes}
\runtitle{AMLE for diffusion processes}

\begin{aug}
\author[A]{\fnms{Jinyuan} \snm{Chang}\thanksref{t1}\ead[label=e1]{changjinyuan1986@pku.edu.cn}}
\and
\author[B]{\fnms{Song Xi} \snm{Chen}\corref{}\thanksref{t1,t2}\ead[label=e2]{songchen@iastate.edu}}
\runauthor{J. Chang and S. X. Chen}
\affiliation{Peking University, and Peking University and Iowa State
University}
\address[A]{Guanghua School of Management\\
Peking University\\
Beijing, 100871\\
China\\
\printead{e1}} %adresu isvedimo komanda gale!
\address[B]{Guanghua School of Management \\
\quad and Center for Statistical Science \\
Peking University\\
Beijing, 100871\\
China\\
and \\
Department of Statistics\\
Iowa State University\\
Ames, Iowa 50011-1210\\
USA\\
\printead{e2}}
\end{aug}

\thankstext{t1}{Supported by Center for Statistical Science
at Peking University and a
National Natural Science Foundation of China Zhongdian Program Grant
(No. 11131002).}

\thankstext{t2}{Supported by NSF Grants DMS-06-04563 and
DMS-07-14978.}

% HISTORY:
\received{\smonth{2} \syear{2011}}
\revised{\smonth{8} \syear{2011}}

% ABSTRACT
%
\begin{abstract}
The transition density of a diffusion process {does not admit an
explicit expression} in general, which prevents the full maximum
likelihood estimation (MLE) based on discretely observed sample paths.
A{\"{\i}}t-Sahalia [\textit{J. Finance} \textbf{54} (1999) 1361--1395;
\textit{Econometrica} \textbf{70} (2002) 223--262] proposed asymptotic
expansions to the transition densities of diffusion processes,
which lead to an approximate maximum likelihood estimation (AMLE) for
parameters. Built on A{\"{\i}}t-Sahalia's [\textit{Econometrica}
\textbf{70} (2002) 223--262; \textit{Ann. Statist.} \textbf{36}
(2008) 906--937] proposal and analysis on the AMLE,
% considered the asymptotic behavior of this estimator.
%He pointed that (i) for fixed sampling interval, there exist a
%sequence about the approximate terms depended on the sample size
%such that the difference between the corresponding AMLEs and the
%full MLE convergence to zero in probability; (ii) for fixed
%approximate term, there exist a sequence about the sampling
%intervals depended on the sample size such that the consistency of
%the corresponding AMLEs hold.
we establish the consistency and convergence rate of the AMLE, which
reveal the roles played by the number of terms used in the
asymptotic density expansions and the sampling interval between
successive observations. We find conditions under which the AMLE
has the same asymptotic distribution as that of the full MLE.
A first order approximation to the Fisher information matrix is
proposed.
\end{abstract}

% KEYWORDS
%
\begin{keyword}[class=AMS]
\kwd[Primary ]{62M05}
\kwd[; secondary ]{62F12}.
\end{keyword}
\begin{keyword}
\kwd{Asymptotic expansion}
\kwd{asymptotic normality}
\kwd{consistency}
\kwd{discrete time observation}
\kwd{maximum likelihood estimation}.
\end{keyword}

\end{frontmatter}

%s1 ###
%s1 #&#
\section{Introduction}\label{sec1}
Continuous-time diffusion processes defined by stochastic
differential equations [\citet{r25}, \citet{r29}, \citet{r31}] are the
basic stochastic modeling tools in the
modern financial theory and applications. Diffusion models are
commonly employed to describe the price dynamics of a financial
asset or a portfolio of assets. An eminent application is in
deriving the price of a derivative contract on an asset or a group
of assets.
%(for instance an index of stocks).
The celebrated Black--Scholes--Merton option pricing formula [\citet{r13},
\citet{r26}]
was obtained by assuming that the
underlying asset followed a geometric Brownian motion such that the
log price process of the underlying asset followed an
Ornstein--Uhlenbeck diffusion process.
% which has a constant diffusion component.
The widely used \citet{r37} and \citet{r15} pricing
formulas for the zero coupon bond
were developed based on two specific mean-reverting diffusion
processes with a~constant or the square root [\citet{r21}]
diffusion functions, respectively.
Other pricing formulas have also
been developed for assets defined by other processes; see \citet{r9} and \citet{r17}. In the
implementations of the aforementioned pricing formula,
the parameters of the diffusion processes which describe
the underlying assets dynamics have to be estimated based on empirical
observations. \citet{r35} gave a comprehensive survey on the financial
applications of continuous-time stochastic models which
were largely the diffusion processes. \citet{r18} provided an
overview on nonparametric estimation for diffusion processes. Other
related works include \citet{r11}, \citet{r38}, \citet{r20}, \citet{r19}, \citet{r27} and
\citet{r8}.

There are several challenges to be faced when estimating parameters of
diffusion processes. One challenge is that despite being
continuous-time models, the
processes are only observed at discrete time points rather than
observed continuously over time. The discrete observations prevent
the use of the relatively straightforward likelihood expressions
[\citet{r30}] available for continuously observed diffusion
processes. Another challenge is that despite the fact that the
diffusion processes
are Markovian, their transition densities from one time point to the
next do not have finite analytic expressions, except for only a few
specific processes. This means that the efficient maximum likelihood
estimation (MLE) cannot be readily implemented for most of these
processes.

In ground-breaking works, A{\"{\i}}t-Sahalia (\citeyear{r1}, \citeyear{r2}) established
series expansions to approximate the
transition densities of univariate diffusion processes. Similar
expansions have been proposed for multivariate processes in
\citet{r3}.
These density approximations, as advocated by A{\"{\i}}t-Sahalia, are then
employed to form approximate likelihood functions, which are maximized
to obtain the approximate maximum likelihood estimators (AMLEs).
A{\"{\i}}t-Sahalia (\citeyear{r2}, \citeyear{r3})
demonstrated that the approximate likelihood converges to the true
likelihood as the number of terms in the series expansions goes to
infinity. He also provided some results on the consistency of the AMLEs.
Numerical evaluations of the transition density approximations as
conducted in \citet{r1}, Stramer and Yan (\citeyear{r33}, \citeyear{r34}) and
others, have shown good performance in the numerical approximation of
the underlying transition densities.
%Statistically speaking,
The approach has opened a very accessible route for obtaining parameter
estimators for diffusion processes, and for estimating other quantities
which are functions of the transition density, as commonly encountered
in finance. % of the underlying asset processes.
Indeed, A{\"{\i}}t-Sahalia and Kimmel (\citeyear{r5}, \citeyear{r6})
demonstrated two such applications in stochastic volatility models and
the affine term structure models, respectively.
\citet{r36} provided some results on the AMLE based on the
one-term expansion for the mean-reverting processes. They revealed
that there was an extra leading order bias term in the AMLE due to
the density approximation.

Although the above-mentioned results on the transition density
approximation and the AMLE had been provided, there are some key
questions that remain to be addressed. One is on the consistency of the
AMLE. While A{\"{\i}}t-Sahalia (\citeyear{r2}, \citeyear{r3})
contained some results on consistency,
there is {more} to be explored. There are two key ingredients in
A\"{\i}t-Sahalia's density approximation. One is $J$, the number of
terms used
in the
approximation, and the other is $\delta$, the length
of the sampling interval between successive observations. In this
paper, we study explicitly the roles played by $J$ and $\delta$ on
the consistency of the AMLE, and quantify their roles on the
convergence rate. Another question is under what conditions on $J$
and~$\delta$, does the AMLE have the same asymptotic distribution as the
full MLE. Here, we consider two regimes: (i) $\delta$ is fixed, and
$J \to\infty$; (ii) $J$ is fixed, but $\delta\to0$, representing
two views of asymptotics.
%In this paper, we establish the consistency of the AMLE and obtain the
%convergence rates under the two asymptotic regimes. % which exhibate
%the roles played by $J$ and $\delta$
%We provide conditions on $J$ and $\delta$ such that the AMLE and the
%full MLE have the same asymptotic normal distribution.
In the case of $\delta\to0$, it is found that $J \geq2$ is
necessary to ensure the AMLE having the same asymptotic normality as
the MLE. Like the transition density, the Fisher information
matrix, the quantity that defines the efficiency of the full MLE, is
unknown analytically; even the underlying transition density is
known. We show in this paper that an approximation to the Fisher
information matrix can be obtained based on the one-term density
approximation.
%
%A specific objective is to reveal the impacts of two key ingridents in
%the approximate MLE.
%
%One is the number of terms (denoted by $J$) used in the density
%approximation; and the other is the length of the %sampling interval
%(denoted by $\delta$) between succesive discrete observations of the
%diffusion process.
%To gain insights on the impacts of $J$ and $\delta$, the bias and
%variance of the AMLE are evaluated under two %settings. One allows $J
%whereas $T$, the total %observation time of the process, goes to
%infinity.
%Our study provides more insights than those reported in Tang and Chen
%(2009) with more generality and more specific %results which reveal
%more specifically on the statistical properties of the AMLE.

The paper is organized as follows. In Section \ref{sec2}, we outline the
transition density approximations of A{\"{\i}}t-Sahalia (\citeyear{r1}, \citeyear{r2}).
Some preliminary analysis is needed for studying the AMLE is presented
in Section \ref{sec3}.
%provide an expansion for the AMLE for a given pair of $J$ and $\delta$.
Section \ref{sec4} establishes the consistency and convergence rates of the
AMLE. Asymptotic normality of the AMLE and its equivalence to the
full MLE are addressed in Section \ref{sec5}. Section \ref{sec6} discusses the
approximation for the Fisher information matrix.
% the main results on the bias and variance of the AMLE. Specific
%examples are presented in Section \ref{sec5}.
Simulation results are reported in Section \ref{sec7}. Technical conditions
and details of proofs are relegated to the \hyperref[app]{Appendix}.

%s2 ###
%s2 #&#
\section{Transition density approximation}\label{sec2}

Consider a univariate diffusion process $(X_t)_{t \ge0}$ defined by
a stochastic differential equation
%
%e2.1 ###
%e2.1 #&#
\begin{equation} \label{eq11}
dX_t=\mu(X_t;\theta)\,dt+\sigma(X_t;\theta)\,dB_t,
\end{equation}
where $\mu$ and $\sigma$ are, respectively, the drift and diffusion
functions and $B_t$ is the standard Brownian motion. Both the drift and
diffusion functions are known except for an unknown parameter
vector $\theta$ taking values in a set
$\Theta\subseteq\mathbb{R}^d$.

Given a sampling interval $\delta> 0$, let $f_X(x | x_0, \delta;
\theta)$ be the transition density of $X_{t +\delta}$ given
$X_t = x_0$ for $(x_0,x) \in{\mathcal{X}\times\mathcal{X}}$, where
${\mathcal{X}}$ is the domain of $X_t$. Despite the parametric forms of
the drift and the diffusion functions that are available
in~(\ref{eq11}), a closed-form expression for $f_X(x | x_0, \delta;
\theta)$ is not generally available for most of the processes. In
most cases, the density is only known to satisfy the Kolmogorov
backward and forward partial differential equations. In
ground-breaking works, A{\"{\i}}t-Sahalia (\citeyear{r1}, \citeyear{r2}) proposed asymptotic
expansions to approximate the transition density.

The approach of A{\"{\i}}t-Sahalia is the following. He first
transformed $X_t$ to a diffusion process with unit diffusion
function by
%
%e2.2 ###
%e2.2 #&#
\begin{equation}\label{eqYt}
Y_t=\gamma(X_t;\theta):=
\int^{X_t}\frac{du}{\sigma(u;\theta)},
\end{equation}
which satisfies
$
dY_t=\mu_Y(Y_t;\theta)\,dt+dB_t$, where
\[
\mu_Y(y;\theta)=\frac{\mu(\gamma^{-1}(y;\theta);\theta)}{\sigma(\gamma
^{-1}(y;\theta);\theta)}-\frac{1}{2}\,\frac{\partial\sigma}{\partial
x}(\gamma^{-1}(y;\theta);\theta).
\]
Let $f_Y(y|y_0,\delta;\theta)$ be the transition density of
$Y_{t + \delta}$ given $Y_t = y_0$. The two density functions are
related according to
%
%e2.3 ###
%e2.3 #&#
\begin{equation}\label{eq2}
f_X(x_t|x_{t-1},\delta;\theta)=\sigma^{-1}(x_t;\theta)\cdot
f_Y(\gamma(x_t;\theta) |\gamma(x_{t-1};\theta),
\delta;\theta).
\end{equation}

To ensure convergence of the expansions, A\"{\i}t-Sahalia standardized
$Y_{t+\delta}$ by $Z_{t+\delta} = \delta^{-1/2} (Y_{t+\delta} -
y_0)$. Let $f_Z(z| y_0, \delta; \theta)$ denote the conditional
density of $Z_{t+\delta}$ given $Z_{t} = 0$, which is related to
$f_Y$ by
\[
f_Z(z| y_0, \delta; \theta) = \delta^{1/2} f_Y(\delta^{1/2} z + y_0 |
y_0,\delta; \theta).
\]
%
%whose conditional density denoted by .
Let $\{ H_j(z)\}_{j=1}^{\infty}$ be the Hermite polynomials
\[
H_j(z) = \phi^{-1}(z) \,{d^j \phi(z) \over d z^j},
\]
which are orthogonal with respect to the standard normal density
$\phi$, namely $\int H_j(z) H_k(z) \phi(x) \,dx = 0$ if $j \ne k$. A
formal Hermite orthogonal series expansion to the density $f_Z(z|
y_0, \delta; \theta)$ is
%
%e2.4 ###
%e2.4 #&#
\begin{equation}\label{eqHermite}
f_Z^{H}(z | y_0, \delta; \theta) =
\phi(z) \sum_{j=0}^{\infty} \eta_j(y_0, \delta; \theta) H_j(z),
\end{equation}
where the coefficients
\begin{eqnarray*}
\eta_j(y_0, \delta; \theta) &=& (j !)^{-1} \int H_j(z) f_Z(z| y_0,
\delta; \theta) \,dz \\
%=& (j !)^{-1} \int H_j\lbrace\delta^{-1/2} (y-y_0) \rbrace f_Y(y|
%y_0, \delta; \theta) dy \\
&=&(j !)^{-1} \mathbb{E}\bigl[ H_j\bigl(\delta^{-1/2} (Y_{t+\delta}-y_0)\bigr) |
Y_t=y_0;\theta\bigr].
\end{eqnarray*}
The last conditional
expectation has no analytic expression in general, although it may
be simulated using the method proposed in \citet{r10}.
A\"{\i}t-Sahalia proposed
Taylor expansions for this conditional expectation with respect to
the sampling interval $\delta$ based on the infinitesimal generator
of $Y_t$. For twice continuously differentiable function $g$, the
infinitesimal generator of $Y_t$ is
%
%e2.5 ###
%e2.5 #&#
\begin{equation}\label{eqgenerator}
\mathcal{A}_{\theta} g(y) = \mu_Y(y; \theta) \,{\partial g \over
\partial y} + \frac{1}{2} \,{\partial^2 g \over\partial
y^2}.
\end{equation}

A $K$-term Taylor series expansion to $\mathbb{E} [
H_j(\delta^{-1/2} (Y_{t+\delta}-y_0)) | Y_t=y_0;\theta]$ is
%
%e2.6 ###
%e2.6 #&#
\begin{eqnarray} \label{eqTaylor}
&&\mathbb{E}\bigl[ H_j\bigl(\delta^{-1/2} (Y_{t+\delta}-y_0)\bigr) | Y_t=y_0;\theta
\bigr]\nonumber\\
&&\qquad= \sum_{k=0}^K \mathcal{A}^k_{\theta} H_j\bigl(\delta^{-1/2}
(y-y_0)\bigr)
\big|_{y=y_0} \frac{\delta^k}{k !}\\
&&\qquad\quad{}+ \mathbb{E}\bigl[\mathcal{A}_{\theta}^{k+1}
H_j\bigl(\delta^{-1/2} (Y_{t+\delta^*}-y_0)\bigr) \big|
Y_t=y_0;\theta\bigr] \frac{\delta^{k+1}}{(k+1) !}. \nonumber
\end{eqnarray}

Substituting (\ref{eqTaylor}) to the orthogonal
expansion (\ref{eqHermite}) followed by gathering terms according
to the powers of $\delta$, a $J$-term expansion to the
transition density $f_Y(y, \delta| y_0; \theta)$ is
\[
f^{(J)}_Y(y |y_0, \delta;\theta) = \delta^{-1/2} \phi\biggl({y-y_0
\over\delta^{1/2}}\biggr) \exp\biggl(\int_{y_0}^y \mu_Y(u; \theta)
\,du \biggr) \sum_{j=0}^J c_{j}(y|y_0; \theta)\frac{ \delta^j}{j !} ,
\]
where $c_0(y|y_0 ; \theta) \equiv1$ and for $j \geq1$,
\begin{eqnarray*}
c_j(y|y_0;\theta)&=&j(y-y_0)^{-j}\\
&&{}\times\int^y_{y_0}(w-y_0)^{j-1}\\
&&\hspace*{26pt}{}\times\biggl\{\lambda_Y(w;\theta)c_{j-1}(w|y_0;\theta)+\frac{1}{2}\,\frac
{\partial^2c_{j-1}(w|y_0;\theta)}{\partial
w^2}\biggr\}\,dw.
\end{eqnarray*}
Here
$\lambda_Y(y;\theta)=-\{\mu^2_Y(y;\theta)+\partial\mu_Y(y;\theta
)/\partial
y\}/2$.

Transforming back from $y$ to $x$ via (\ref{eqYt}) and
(\ref{eq2}), the $J$-term expansion to $f_X(x|x_0, \delta; \theta)$
is
\begin{eqnarray*}
&&
f^{(J)}_X(x|x_0,\delta; \theta) \\
&&\qquad= \sigma^{-1}(x; \theta) \delta^{-1/2}
\phi\biggl({\gamma(x; \theta)- \gamma(x_0; \theta)\over\delta^{1/2}}\biggr) \\
&&\qquad\quad{} \times \exp\biggl\{\int_{x_0}^x
\frac{\mu_Y(\gamma(u;\theta);\theta)} {\sigma(u;\theta)} \,du \biggr\}
\sum_{j=0}^J c_{j}(\gamma(x; \theta) |\gamma(x_0; \theta); \theta)
\frac{\delta^j}{j!} . %\label{eqJ-exp}
\end{eqnarray*}
Although it employs the Hermite polynomials and has the Gaussian
density as the leading term as an Edgeworth expansion does, the
transition density expansion is not an Edgeworth expansion. This
is because the latter is for density functions of statistics
admitting the central limit theorem, which differs from the current
context of expanding the transition density.
\citet{r2} demonstrated that
as $J \to\infty$,
%
%e2.7 ###
%e2.7 #&#
\begin{equation} \label{eqconv1}
f_X^{(J)}(x | x_0, \delta; \theta) \to f_X(x
| x_0, \delta; \theta)
\end{equation}
uniformly with respect
to $\theta\in\Theta$ and $x_0$ over compact subsets of
${\mathcal{X}}$. The convergence is also uniform with respect to $x$ over
subsets of ${\mathcal{X}}$ depending on the property of $\sigma(x;
\theta)$.

%under quite conditions that require the smoothness of the drift and
%diffusion functions, the nondegeneracy of $\sigma^2(\cdot; \theta)$
%and %condtions that specify the growth rates on $\mu_Y(y;\theta)$ near
%the boundaries together with a bounded potential function.
%the that this $J$-term expansion converges to $f_x(x, \delta|x_0;
Define
\begin{eqnarray*}
A_1(x|x_0, \delta;\theta)&=&-\log
\{\sigma(x;\theta)\}-\frac{1}{2\delta}
\{\gamma(x;\theta)-\gamma(x_0;\theta)\}^2,\\
A_2(x|x_0, \delta;\theta) &=&
\int_{x_0}^{x}\frac{\mu_Y(\gamma(u;\theta);\theta)}
{\sigma(u;\theta)}\,du
\end{eqnarray*}
and
\[
A_3(x|x_0, \delta;\theta) =\log
\Biggl\{\sum_{j=0}^{J}c_j(\gamma(x; \theta)| \gamma(x_0;
\theta);\theta)\delta^j/j!\Biggr\}.
\]
%
%then
%&\log f_X^{(J)}(x|x_0, \delta;\theta) \\
%&=-\log\sqrt{2\pi\delta}+A_1(x|x_0,
If ${\sum_{j=0}^\infty}
|c_j(y|y_0,\delta;\theta)|\delta^j/j!<\infty$ on $\mathcal
{Y}\times\mathcal{Y}$ with probability one, where $\mathcal{Y}$ is
the domain of $Y_t$, we can define $\tilde{A}_3(x|x_0,
\delta;\theta)=\log
\{\sum_{j=0}^{\infty}c_j(y|y_0;\theta)\delta^j/j!\}$. Then the
result in (\ref{eqconv1}) implies that
%
%e2.8 ###
%e2.8 #&#
\begin{eqnarray}\label{eqlog-exp}
&&\log f_X(x|x_0, \delta;\theta) \nonumber\\
&&\qquad=-\log\sqrt{2\pi\delta}+A_1(x|x_0,
\delta;\theta)+A_2(x|x_0, \delta;\theta)\\
&&\qquad\quad{}+\tilde{A}_3(x|x_0,
\delta;\theta).\nonumber
\end{eqnarray}
Expression (\ref{eqlog-exp}) is the starting point for our analysis.
%%In
%Appendix, we will give some conditions to guarantee
%$\tilde{A}_3(x|x_0,\delta;\theta)$ has some properties which is
%very important in our research.

Given a set of discrete observations $\{X_{t \delta}\}_{t=1}^n$ with
equal sampling length~$\delta$ of the diffusion process $(X_t)_{t
\ge0}$, to simplify notations, we write $X_t$ for~$X_{t \delta}$
and hide $\delta$ in the expressions for the transition density
$f_X$ and its approximations. At the same time, we use $f$ and
$f^{(J)}$ to express $f_X$ and $f_X^{(J)}$, respectively. Based on
the $J$-term expansion to the true transition density, the $J$-term
approximate log-likelihood function given in \citet{r2} is
\begin{eqnarray*}
{\ell_{n,\delta}^{(J)}}(\theta)%=&\sum_{t=1}^n \log f^{(J)}(X_{t} |
%X_{t-1},\delta; \theta)\\
&=&- n \log\sqrt{2\pi\delta}+ \sum_{t=1}^n
A_1(X_t|X_{t-1},\delta;\theta)\\
&&{} + \sum_{t=1}^n A_2(X_t|X_{t-1},\delta;\theta)+ \sum_{t=1}^n
A_3(X_t|X_{t-1},\delta;\theta).
\end{eqnarray*}

Let ${\hat{\theta}_{n,\delta}^{(J)}} = \arg\max_{\theta\in\Theta}
{\ell_{n,\delta}^{(J)}}(\theta)$ be the approximate MLE (AMLE) and
${\hat{\theta}_{n,\delta}}$ be the true MLE that maximizes the full
likelihood
\[
{\ell_{n,\delta}}(\theta)=\sum_{t=1}^n \log f(X_{t} |
X_{t-1},\delta; \theta).
\]
To keep the notation simple, we write
$\hat{\theta}^{(J)}_n = \hat{\theta}_{n,\delta}^{(J)}$ and
$\hat{\theta}_n = \hat{\theta}_{n,\delta}$ by suppressing $\delta$
in subscripts.

%consistency are the same in the two cases. I plan to delete this
%paragraph.

%s3 ###
%s3 #&#
\section{Preliminaries}\label{sec3}

Under regular circumstances as assumed by condi-\break tion~(A.2)(ii)~in~the
\hyperref[app]{Appendix}, the full MLE $\hat{\theta}_{n}$ and the $J$-term
approximate MLE $\hat{\theta}^{(J)}_{n}$ satisfy their respective
likelihood score equations so that
%
%e3.1 ###
%e3.1 #&#
\begin{equation}\label{eq03}
\sum_{t=1}^n\nabla_\theta\log
f(X_t|X_{t-1},\delta;{\hat{\theta}_{n}})=\sum_{t=1}^n\nabla_\theta\log
f^{(J)}\bigl(X_t|X_{t-1},\delta;{\hat{\theta}^
{(J)}_{n}}\bigr)=0.\hspace*{-22pt}
\end{equation}
Subtracting $\sum_{t=1}^n \nabla_\theta\log
f^{(J)}(X_t|X_{t-1},\delta;\theta_0)$ from both sides of
(\ref{eq03}),
%
%e3.2 ###
%e3.2 #&#
\begin{eqnarray}\label{eq04}
&&\sum_{t=1}^n\nabla_\theta\log
f^{(J)}\bigl(X_t|X_{t-1},\delta;{\hat{\theta}^
{(J)}_{n}}\bigr)-\sum_{t=1}^n\nabla_\theta\log
f^{(J)}(X_t|X_{t-1},\delta;\theta_0)\nonumber\\
&&\qquad=\sum_{t=1}^n\nabla_\theta[\tilde{A}_3(X_t|X_{t-1},\delta;\theta
_0)-A_3(X_t|X_{t-1},\delta;\theta_0)]
\\
&&\qquad\quad{}+\sum_{t=1}^n\nabla_\theta\log
f(X_t|X_{t-1},\delta;{\hat{\theta}_{n}})-\sum_{t=1}^n\nabla_\theta\log
f(X_t|X_{t-1};\theta_0).
\nonumber
\end{eqnarray}
Carrying out Taylor expansions on both sides of (\ref{eq04}), we
can get
%
%e3.3 ###
%e3.3 #&#
\begin{eqnarray}\label{eqtaylor-eq}\quad
&&\frac{1}{n} \sum_{t=1}^n\nabla^2_{\theta\theta}\log
f^{(J)}(X_t|X_{t-1},\delta;\theta_0)\cdot\bigl({\hat{\theta
}^{(J)}_{n}}-\theta_0\bigr) \nonumber\\
&&\quad{}+\frac{1}{2} \bigl[E_d\otimes\bigl({\hat{\theta}_{n}^{(J)}}-\theta_0\bigr)'\bigr]\cdot
\frac{1}{n}\sum_{t=1}^n\nabla^3_{\theta\theta\theta}\log
f^{(J)}(X_t|X_{t-1},\delta;\tilde{\theta})\cdot\bigl({\hat{\theta
}_{n}^{(J)}}-\theta_0\bigr) \nonumber\\
&&\qquad=\frac{1}{n} \sum_{t=1}^{n}\nabla_\theta[\tilde
{A}_3(X_t|X_{t-1},\delta;\theta_0)-A_3(X_t|X_{t-1},\delta;\theta_0)]\\
&&\qquad\quad{}+\frac{1}{n}\sum_{t=1}^{n}\nabla^2_{\theta\theta}\log
f(X_t|X_{t-1},\delta;\theta_0) \cdot({\hat{\theta}_{n}}-\theta_0) \nonumber\\
&&\qquad\quad{}+\frac{1}{2} [E_d\otimes
({\hat{\theta}_{n}}-\theta_0)']\cdot
\frac{1}{n}\sum_{t=1}^n\nabla^3_{\theta\theta\theta}\log
f(X_t|X_{t-1},\delta;\bar{\theta})\cdot({\hat{\theta}_{n}}-\theta_0),
\nonumber
\end{eqnarray}
where $E_d$ is the $d\times d$ identity matrix, $\tilde{\theta}$ is
on the joint line between ${\hat{\theta}_{n}^{(J)}}$ and~$\theta_0$
and $\bar{\theta}$ is on the joint line between $\hat{\theta}_{n}$
and $\theta_0$. Here we define
\[
\nabla^3_{\theta\theta\theta}\log
f(X_t|X_{t-1},\delta;\theta):=\pmatrix{
\partial^3\log
f(X_t|X_{t-1},\delta;\theta)/\partial\theta\,\partial\theta'\,\partial\theta
_1 \cr
\vdots\cr
\partial^3\log
f(X_t|X_{t-1},\delta;\theta)/\partial\theta\,\partial\theta'\,\partial\theta
_d },
\]
which is a $d^2\times d$ matrix, and
$\nabla^3_{\theta\theta\theta}\log
f^{(J)}(X_t|X_{t-1},\delta;\theta)$ is similarly defined.
Furthermore, let
\begin{eqnarray*}
F_n(\theta_0,J,\delta) &=& n^{-1}
\sum_{t=1}^{n}\nabla^2_{\theta\theta}
[\tilde{A}_3(X_t|X_{t-1},\delta;\theta_0)-A_3(X_t|X_{t-1},\delta
;\theta_0)], \\
U_n(\theta_0,J,\delta) &=& n^{-1}
\sum_{t=1}^{n}\nabla_{\theta}
[\tilde{A}_3(X_t|X_{t-1},\delta;\theta_0)-A_3(X_t|X_{t-1},\delta
;\theta_0)]
\end{eqnarray*}
and
\[
N_n(\theta_0,J,\delta)= n^{-1}
\sum_{t=1}^{n}\nabla^2_{\theta\theta}\log
f^{(J)}(X_t|X_{t-1},\delta;\theta_0).
\]
Then (\ref{eqtaylor-eq}) can be written as
%
%e3.4 ###
%e3.4 #&#
\begin{eqnarray}\label{eqkeyeq}
&&N_n(\theta_0,J,\delta)\bigl({\hat{\theta}_{n}^{(J)}}-\theta_0\bigr)+ \Delta_{n
1}\bigl({\hat{\theta}^{(J)}_{n}}, \theta_0\bigr) \nonumber\\
&&\qquad=U_n(\theta_0,J,\delta) + [N_n(\theta_0,J,\delta) +F_n(\theta
_0,J,\delta) ]({\hat{\theta}_{n}}-\theta_0)\\
&&\qquad\quad{}+ \Delta_{n 2}({\hat{\theta
}_{n}}, \theta_0),
\nonumber
\end{eqnarray}
where
$\Delta_{n 1}({\hat{\theta}^{(J)}_{n}}, \theta_0)$
and $\Delta_{n 2}({\hat{\theta}_{n}}, \theta_0)$ denote the
{remainder terms whose explicit expressions can be obtained by
matching (\ref{eqtaylor-eq}) with (\ref{eqkeyeq}).
% Then,
%(\ref{eqtaylor-eq}) can be written as

Expansion (\ref{eqkeyeq}) is the starting point in our studies
for the consistency and asymptotic distribution of the AMLE.
%$\bf\hat{\theta}_{n}^{(J)}$ for (i) fixed $\delta$ but $J \to
% and $ obtaining the central limit
%theorem for $\hat{\theta}_n^{(J)}$ when $J$ is fixed but $J \ge2$.
%A more elaborate expansion than (\ref{eqkeyeq}) with quadratic
%terms will be provided for the case of $J=1$ in the next section.
%%while $\delta\to0$.}}
Indeed, the asymptotic properties of the AMLE will be evaluated under
two regimes regarding $J$ and $\delta$.
The first one is that
%
%e3.5 ###
%e3.5 #&#
\begin{equation}\label{eqregime1}
\mbox{$\delta$ is fixed} \qquad\mbox{but $J \to\infty$,}
\end{equation}
which is the situation considered in \citet{r2}.
The second regime allows that
%
%e3.6 ###
%e3.6 #&#
\begin{equation} \label{eqregime2}
\mbox{$J$ is fixed},\qquad \mbox{$\delta\rightarrow0$}\qquad\mbox{but $n\delta\rightarrow\infty
$,}
\end{equation}
which is more tuned with an implementation of the density
approximation with a fixed number of terms.

We will first present some results which are valid for any fixed $J$
and $\delta$.
%We firstly introduce a matrix norm
Let $\|A\|_2=\{\rho(A'A)\}^{1/2}$ be the spectral norm of a matrix
$A$, where $\rho(A'A)$ denotes the largest eigen-value of $A' A$.
The following proposition describes properties for the quantities that
appear in (\ref{eqkeyeq}).
%
%pr1 #&#
\begin{pn}\label{prop1}
Under conditions \textup{(A.1), (A.3), (A.4), (A.6), (A.7)} given in the \hyperref[app]{Appendix},
there exists a positive constant $\Delta$ such that for any positive
integer $J$ and $\delta\in(0,\Delta)$:\vspace*{8pt}

\textup{(a)} $\mathbb{E}\{F_n(\theta_0,J,\delta)\}$,
$\mathbb{E}\{U_n(\theta_0,J,\delta)\}$ and
$\mathbb{E}\{N_n(\theta_0,J,\delta)\}$ exist;\vspace*{1pt}

\textup{(b)} $ \Delta_{n 1}({\hat{\theta}^{(J)}_{n}},
\theta_0)=O_p\{\|{\hat{\theta}_{n}^{(J)}}-\theta_0\|_2^2\}$ and
$\Delta_{n 2}({\hat{\theta}_{n}},
\theta_0)=O_p\{\|{\hat{\theta}_{n}}-\theta_0\|_2^2\}$ as $n \to
\infty$.
\end{pn}

Let $I(\delta) = -\mathbb{E}\nabla^2_{\theta\theta}\log
f(X_t|X_{t-1},\delta;\theta_0)$ be the Fisher information matrix,
which we assume is invertible in condition (A.5). It is expected
that the expected value of $N_n(\theta_0, J, \delta)$, denoted by
$N(\theta_0, J, \delta)$, will converge to $-I(\delta)$, as $J \to
\infty$ for each fixed $\delta$ or $J$ being fixed but
$\delta\rightarrow0$. The following proposition bounds the
difference between $N(\theta_0, J, \delta)$ and $-I(\delta)$ for
each fixed $J$ and $\delta$.
%
%pr2 #&#
\begin{pn}\label{prop2}
Under conditions \textup{(A.1), (A.4), (A.6), (A.7)} given in the
\hyperref[app]{Appendix}, there exist two positive constants $\bar{\Delta}$ and $C$,
that are not dependent on $J$ and $\delta$,
such that for any positive integer $J$ and
$\delta\in(0,\bar{\Delta})$,
\[
\|N(\theta_0,J,\delta)+I(\delta)\|_2\leq C\delta^{J+1}.
\]
\end{pn}

As $I(\delta)$ is invertible for each fixed
$\delta
>0$, $N_n(\theta_0, J, \delta)$ will be invertible with probability
approaching one as $J \to\infty$ for a fixed $\delta$. However,
{if} $\delta\to0$, the limit of the Fisher information
$I(0):=\lim_{\delta\to0}I(\delta)$, as well as~$N(\theta_0, J,
0)$, may be singular. This is the case for some Ornstein--Uhlenbeck
processes as shown in Section~\ref{sec6}. The following proposition provides
another account on $N(\theta_0,J,\delta)$ and its deviation from
$-I(\delta)$, as well as the convergence of
$N^{-1}(\theta_0,J,\delta)U(\theta_0,J,\delta)$, where
$U(\theta_0,J,\delta)$ denotes the expected value of~$U_n(\theta_0,J,\delta)$ for each pair of fixed $J$ and $\delta$.
%
%pr3 #&#
\begin{pn}\label{prop3}
Under conditions \textup{(A.1), (A.3)--(A.7)} given in the \hyperref[app]{Ap-}
\hyperref[app]{pendix}, there exist
two constants $C_1, C_2$, that are not dependent on $J$ and
$\delta$, and a constant $\underline{\Delta}>0$ such that for any
positive integer $J$ and $\delta\in(0,\underline{\Delta})$,
\[
\|N^{-1}(\theta_0,J,\delta)I(\delta)+E_d\|_2\leq
C_1\delta^J \quad\mbox{and}\quad\|N^{-1}(\theta_0,J,\delta)U(\theta_0,J,\delta
)\|_2\leq
C_2\delta^J.
\]
%
%under either (i) for any fixed $\delta\in(0,\bar{\Delta})$, where
%$\bar{\Delta}$ is the quantity in Proposition 2, and
%$J\rightarrow\infty$ or (ii) for any fixed $J$, $\delta
\end{pn}

The next proposition describes the convergence rate for the
difference between the first derivatives of the full log-likelihood
and the approximate log-likelihood.
%
%pr4 #&#
\begin{pn}\label{prop4}
Under\vspace*{1pt} conditions \textup{(A.1), (A.4), (A.6), (A.7)} given in the \hyperref[app]{Appendix}, there
exist two finite\vspace*{1pt} positive constants $\tilde{\Delta}$ and $C$,
not dependent on $J$ and~$\delta$, such that for any $J$,
$\delta\in(0,\tilde{\Delta}]$ and $n$,
\[
\mathbb{E}\Bigl\{\sup_{\theta\in\Theta}\bigl\|n^{-1}\cdot\nabla_\theta\bigl[\ell
_{n,\delta}(\theta)-\ell_{n,\delta}^{(J)}(\theta)\bigr]\bigr\|_2\Bigr\}\leq
C\delta^{J+1}.
\]
\end{pn}

The following proposition together with Proposition \ref{prop4} %\footnote{I
%suggest we move these two prop. to the appendix.}
is needed to establish the consistency of the AMLE.
%
%pr5 #&#
\begin{pn}\label{prop5}
Under conditions \textup{(A.1), (A.3), (A.4), (A.6), (A.7)} given in the \hyperref[app]{Appendix},
there exists a constant $\dot{\Delta}>0$ such that
\[
\sup_{\theta\in\Theta}\Biggl\|\frac{1}{n}\sum_{t=1}^n\nabla_\theta\log
f(X_t|X_{t-1},\delta;\theta)-\mathbb{E}\nabla_\theta\log
f(X_t|X_{t-1},\delta;\theta)\Biggr\|_2\stackrel{p}{\rightarrow}0
\]
for \textup{(i)} $\delta\in(0,\dot{\Delta}]$ being fixed,
$n\rightarrow\infty$, or \textup{(ii)} $n\rightarrow\infty$,
$\delta\rightarrow0$ but $n\delta\rightarrow\infty$.
\end{pn}

As the full MLE $\hat{\theta}_{n}$ is a key bridge for the AMLE, we
report in the following proposition
the asymptotic normality of the
MLE which covers both cases of fixed $\delta$ and diminishing
$\delta$ case. %The result is not new. Our reporting here is just for
%easy referencing.
%
%pr6 #&#
\begin{pn}\label{prop6}
Under conditions \textup{(A.1)--(A.7)}
given in the \hyperref[app]{Appendix},
\[
\sqrt{n}I^{1/2}(\delta)(\hat{\theta}_n-\theta_0)\stackrel{d}{\rightarrow
}N(0,E_d)\qquad\mbox{as }n\delta^3\rightarrow\infty,
\]
where $E_d$ is $d\times d$ identity matrix.
\end{pn}

The requirement of $n\delta^3 \to\infty$ in the above proposition
is to cover the case where $I(0) = \lim_{\delta\to0} I(\delta)$ is
singular, as spelled out in the proof given in the \hyperref[app]{Appendix}. If such
case is ruled out, for instance, via the so-call Jacobsen condition
[\citet{r24}, \citet{r32}], the more standard $n \delta
\to\infty$ is sufficient; see also \citet{r23} for related
results.

%considered a general theory of efficient
%estimation for ergodic diffusions sampled at high frequency. In this
%paper, he got the estimation based on a general estimating function
%which is exactly or approximately a martingale estimating function.
%This idea includes the full MLE. Under some conditions mentioned in
%this paper, we can find the rates of convergence for estimators of
%parameters in the drift coefficients are $(n\delta)^{-1/2}$, whereas
%that for estimators of parameters in the diffusion coefficients are
%$n^{-1/2}$. He called these rates are optimal rates. But the
%asymptotic variance mentioned in his paper, if we choose the score
%function to be the estimating function, is not the Fisher
%information matrix. Above proposition reveals under what conditions,
%the full MLE will obtain the lower bound of asymptotic variance. }

%There are three parts in this section. Firstly, we consider Secondly,
%we study the central
%limit theorem of $\hat{\theta}_n^{(J)}$. At
%last, we consider the leading order of bias and variance of
%$\hat{\theta}_n^{(J)}$.

%s4 ###
%s4 #&#
\section{Consistency}\label{sec4}
We consider in this section the consistency of
the\break AMLE~$\hat{\theta}_{n}^{(J)}$ and establish its convergence rate under
the two asymptotic regimes given in (\ref{eqregime1}) and
(\ref{eqregime2}), respectively. The two asymptotic regimes were
also considered in A\"{\i}t-Sahalia (\citeyear{r2}, \citeyear{r3}). For a fixed
sampling interval~$\delta$, \citet{r2}  proved that
there existed a sequence $J_n \to\infty$ such that
${\hat{\theta}_{n}^{(J_n)} - \hat{\theta}_{n}} \stackrel{p} \to0$
under $P_{\theta_0}$ as $n \to\infty$, where $P_{\theta_0}$ is the
underlying probability measure.
%consistency; (ii) What does he mean by ``under $P_{\theta_0}$" ?
%Should it be the %conditional prob measure ?}.
Based on the
consistency of $\hat{\theta}_{n}$, we know that the consistency of
$\hat{\theta}_{n}^{(J_n)}$ is hold.
%But how to find relationship between $J_n$
%and $n$ wasn't mentioned in this paper.
For a fixed $J$, \citet{r3} proved that there existed a
sequence $\{\delta_n\}$ vanishing to zero such that
$\sqrt{n}I^{1/2}(\delta_n)({\hat{\theta}_{n,\delta_n}^{(J)}}-\theta_0)=O_p(1)$.

In this paper, we will give more explicit guidelines on how to
select the afore-mentioned sequences $J_n$ and $\delta_n$ so that
the AMLE is consistent. Our study here begins with (\ref{eq03}),
which together with Propositions \ref{prop4} and \ref{prop5} lead to the following {result}
on the consistency of the AMLE under the two asymptotic
regimes, respectively.
%
%th1 #&#
\begin{tm}\label{theo1}
Under conditions \textup{(A.1)--(A.4), (A.6), (A.7)} given in the \hyperref[app]{Appendix},
$\hat{\theta}_{n}^{(J)}-\theta_0\stackrel{p}{\rightarrow}0$ under
either: \textup{(i)}
$\delta\in(0,\tilde{\Delta}\wedge\dot{\Delta}]$ being fixed,
$J\rightarrow\infty$ and $n\rightarrow\infty$, or \textup{(ii)} $J$ being
fixed, $n\rightarrow\infty$, $\delta\rightarrow0$ but
$n\delta\rightarrow\infty$.\vadjust{\goodbreak}
\end{tm}

By Proposition \ref{prop2} and condition (A.5),
% we know that the matrix $N(\theta_0,J,\delta)$ is invertable when $J$
%is sufficient large
%or $\delta$ is sufficient small.
multiply $N^{-1}(\theta_0,J,\delta)$ on both sides of~(\ref{eqkeyeq}),
we have
%
%e4.1 ###
%e4.1 #&#
\begin{eqnarray}\label{eq41}
&&{\hat{\theta}_{n}^{(J)}}-\theta_0\nonumber\\
&&\qquad=N^{-1}U_n+N^{-1}(N_n+F_n)({\hat{\theta}_{n}}-\theta
_0)-N^{-1}(N_n-N)\bigl({\hat{\theta}_{n}^{(J)}}-\theta_0\bigr)\\
&&\qquad\quad{}-N^{-1}\Delta_{n1}\bigl({\hat{\theta}_{n}^{(J)}},\theta_0\bigr)+N^{-1}\Delta
_{n2}({\hat{\theta}_{n}},\theta_0).\nonumber
\end{eqnarray}
%
%This equation is the beginning of our following work.
From this together with Proposition \ref{prop4} and Theorem \ref{theo1}, we can
establish the convergence rate of the AMLE.
%
%th2 #&#
\begin{tm}\label{theo2}
Under conditions \textup{(A.1)--(A.7)} given
in the \hyperref[app]{Appendix},
\[
{\hat{\theta}_{n}^{(J)}}-\theta_0=\cases{
O_p\{\delta^{J+1}+(n\delta)^{-1/2}\}, &\quad if $\delta\in
(0,\tilde{\Delta}\wedge\dot{\Delta}]$ is fixed and $J \to\infty$;
\vspace*{2pt}\cr
O_p\{\delta^{J}+(n\delta)^{-1/2}\}, &\quad if $J$ is fixed, $\delta
\rightarrow0$ but $n\delta^3\rightarrow\infty$.}
\]
\end{tm}

The above theorem reveals the impacts of the sampling interval
$\delta$ and the number of terms $J$ used in the density
approximation on the convergence rate. In particular, the rate of
AMLE has an extra $\delta^{J+1}$ or $\delta^J$ term in addition to
the standard rate $(n \delta)^{-1/2}$ of the full MLE. This extra
term is the result of the density approximation, and its particular
form suggests\vspace*{1pt} that the sampling interval $\delta$ has to be less
than 1 in order to make the AMLE~$\hat{\theta}_{n}^{(J)}$ converge
to~$\theta_0$. It is apparent\vspace*{1pt} that the higher the $J$ is, the less
impact the extra term has on the AMLE $\hat{\theta}_{n}^{(J)}$.

%s5 ###
%s5 #&#
\section{Asymptotic distribution}\label{sec5}

In this section, we consider the asymptotic distribution of the AMLE
$\hat{\theta}_{n}^{(J)}$. Our investigations are organized according
to two asymptotic regimes: (i) $\delta$ fixed, $J\rightarrow\infty$
and (ii) $J$ fixed, $\delta\rightarrow0$ but
$n\delta\rightarrow\infty$.

%s5.1 ###
%s5.1 #&#
\subsection{\texorpdfstring{Fixed $\delta$, $J \to\infty$}{Fixed delta, J to infinity}}\label{sec51}

This is a simple case to treat. Under this setting, we note from
Proposition \ref{prop2} and condition (A.5) that
$N^{-1}(\theta_0,J,\delta)=O(1)$ uniformly for any $J$. Utilizing
the result in Theorem \ref{theo2}, expansion (\ref{eq41}) becomes
\[
{\hat{\theta}_{n}^{(J)}}-\theta_0
%=&N^{-1}U_n+N^{-1}(N_n+F_n)({\bf\hat{\theta}_{n}}-
%=&N^{-1}U_n-N^{-1}I({\bf\hat{\theta}_{n}}-\theta_0)+O_p(n^{-1/2}
=N^{-1}U_n+({\hat{\theta}_{n}}-\theta_0)+O_p(n^{-1/2}\delta
^{J-1/2}+n^{-1}\delta^{-1} + \delta^{2J+2}).
\]
%
%We recognize that the first term on the right hand side is associated
%with $U_n = O_p(\delta^{J+1})$, as a result of the approximation as
%well as having $N(\theta_0,J,\delta)$ based on the $J$-term instead of
%$-I(\delta)$.
Hence, note that $U_n=O_p(\delta^{J+1})$,
\begin{eqnarray*}
&&\sqrt{n}I^{1/2}(\delta)
\bigl({\hat{\theta}_{n}^{(J)}}-\theta_0\bigr)\\
&&\qquad=\sqrt{n}I^{1/2}(\delta) ({\hat{\theta}_{n}}-\theta_0)+O_p(\delta
^{J-1/2}+n^{-1/2}\delta^{-1} + n^{1/2}\delta^{J+1}).
%=&\sqrt{n}I^{1/2}(\delta)(\hat{\theta}_n-\theta_0)+O_p(n^{1/2}
\end{eqnarray*}
If $n\delta^{2J+2}\rightarrow0$, then
\[
\sqrt{n}I^{1/2}(\delta)\bigl({\hat{\theta}_{n}^{(J)}}-\theta_0\bigr)\stackrel
{d}{\rightarrow}N(0,E_d).
\]
Therefore, the AMLE has the same asymptotic distribution as the full
MLE~$\hat{\theta}_{n}$.
This is attained by requesting $n\delta^{2J+2}\rightarrow0$ in addition
to $J \to\infty$. %This condition implicitly imply $\delta< 1$.
If $n \delta^{2J+2} \to c > 0$, the AMLE is still
asymptotically normal but would have an inflated variance due to the
contribution from the first term involving $U_n$. Apart from this,
the asymptotic mean will no longer be zero. Hence, it is much desirable
to have $ n \delta^{2J+2}
\to0$. The latter condition prescribes a rule on the selection
of the $J=J_n(\delta)$. By choosing an $\varepsilon> 0$ so that
$\delta^{2J+2}=n^{-1-\varepsilon}$ for each pair of $n$ and $\delta$,
then
\[
J=J_n(\delta)= \frac{-1-\varepsilon}{2\log\delta}\log
n-1>\frac{-1}{2\log\delta}\log n-1.
\]
The integer truncation of the above lower bound plus one can be used
as a~reference value for the number of terms used in the density
approximation for each given pair of $(n,\delta)$.

%t1 ###
%t1 #&#
\begin{table}
\tablewidth=240pt
\caption{The least approximation term selection to guarantee the
AMLE has the same asymptotic distribution as the full MLE for
special sampling interval $\delta$ and sample size $n$}
\label{table1}
\begin{tabular*}{\tablewidth}{@{\extracolsep{\fill}}lcccc@{}}
\hline
% after \ \hline or \cline{col1-col2} \cline{col3-col4} ...
$\bolds{\delta}$ & $\bolds{n=500}$ & $\bolds{n=1\mbox{\textbf{,}}000}$
& $\bolds{n=2\mbox{\textbf{,}}000}$ & $\bolds{n=4\mbox{\textbf{,}}000}$\\
\hline
$1/252$ & \hphantom{0}1 & \hphantom{0}1 & \hphantom{0}1 &\hphantom{0}1\\
$1/52$ & \hphantom{0}1 & \hphantom{0}1 & \hphantom{0}1&\hphantom{0}1\\
$1/12$ & \hphantom{0}1 & \hphantom{0}1 & \hphantom{0}1&\hphantom{0}1\\
$1/4$ & \hphantom{0}2 & \hphantom{0}2 & \hphantom{0}2 &\hphantom{0}2\\
$1/2$ & \hphantom{0}4 & \hphantom{0}4 & \hphantom{0}5 &\hphantom{0}5\\
$3/4$ & 10 & 12 & 13 &14\\
\hline
\end{tabular*}
\end{table}

Table \ref{table1} reports such reference values of $J$ assigned by the above
formula for a set of $(n,\delta)$ combinations commonly encountered
in empirical studies. It shows that for monthly frequency or less
($\delta\leq1/12$), one term approximation is adequate, and for
$\delta=1/4$, $J=2$ is needed. However, there is a dramatic increase
in $J$ as the sampling length is larger than $1/4$: demanding at
least four terms for $\delta=1/2$ (half yearly) or at least ten
terms for $\delta= 3/4$. The number of terms also increases for
these higher $\delta$ values as $n$ increases, although the rate of
this increase is much slower than that as $\delta$ is increased. The
latter may be understood that for a given~$\delta$, as $n$
increases, the chance of having extreme values in the tails of the
transition distribution increases. As the density approximation
is less accurate in the tails than in the main body of the
distribution, there is a need for having more terms in the density
approximation.

%Above inequality gives the selection for the approximation order $J$
% for a fixed $\delta$.

%s5.2 ###
%s5.2 #&#
\subsection{\texorpdfstring{$J$ fixed, $\delta\rightarrow0$ but $n\delta\rightarrow\infty$}
{J fixed, delta rightarrow 0 but n delta rightarrow infinity}}
\label{sec52}

Our starting point is the expansion~(\ref{eq41}). As
$N_n-N=O_p\{(n\delta)^{-1/2}\}$, $N^{-1} (N_n - N) = o_p(1)$ if $n
\delta^3 \to\infty$, which is also required in the asymptotic
normality of the full MLE
as outlined in Proposition \ref{prop6}.
%I make a mistake in
%the original version, there should be $n\delta^3\rightarrow\infty$.
%} which we will assume in the
%rest of this section when we are dealing with diminishing $\delta$.
{We will show in the following that $n \delta^3 \to\infty$ is also necessary
to ensure AMLE sharing the same asymptotic distribution
as the full MLE. It is understood that in order for
$\hat{\theta}_{n}^{(J)}$ having the same asymptotic distribution as
$\hat{\theta}_{n}$,
%and the definition of
%$\Delta_{n1}(\hat{\theta}_n^{(J)},\theta_0)$ and
%$\Delta_{n2}(\hat{\theta}_n,\theta_0)$, %we know that, e need that
it is required that
\[
N^{-1} U_n,
N^{-1}\Delta_{n1}\bigl({\hat{\theta}_{n}^{(J)}},\theta_0\bigr)\mbox{ and }
N^{-1}\Delta_{n2}({\hat{\theta}_{n}},\theta_0)\qquad\mbox{are all }
o_p\bigl\{\bigl\|{\hat{\theta}_{n}^{(J)}}-\theta_0\bigr\|_2\bigr\}.
\]

We will demonstrate in the following that the above requirements can
be attained by $n \delta^3 \to\infty$ and $J \geq2$. Hence,
under these circumstances, $\hat{\theta}_{n}^{(J)}$ has the same
asymptotic distribution as $\hat{\theta}_{n}$. Later we will
demonstrate that this equivalence in the asymptotic distribution is
quite unlikely for $J=1$.
Our analysis needs to expand (\ref{eq04}) to the quadratic terms.
To this end, let us define
\[
M_n(\theta_0,J,\delta)=n^{-1}\sum_{t=1}^n\nabla^3_{\theta\theta\theta
}\log
f^{(J)}(X_t|X_{t-1},\delta;\theta_0)
\]
\mbox{and}
\[
T_n(\theta_0,J,\delta)=n^{-1}\sum_{t=1}^n\nabla^3_{\theta\theta\theta
}\log
f(X_t|X_{t-1},\delta;\theta_0).
\]
By further expanding to quadratic terms, (\ref{eq41}) can be
written as
%
%e5.1 ###
%e5.1 #&#
\begin{eqnarray}\label{eqtaylor2}\quad
&&{\hat{\theta}_{n}^{(J)}}-\theta_0\nonumber\\
&&\qquad=N^{-1}U_n+N^{-1}(N_n+F_n)({\hat{\theta}_{n}}-\theta
_0)-N^{-1}(N_n-N)\bigl({\hat{\theta}_{n}^{(J)}}-\theta_0\bigr)\nonumber\\
&&\qquad\quad{}-\tfrac{1}{2} N^{-1}\bigl[E_d\otimes\bigl({\hat{\theta}_{n}^{(J)}}-\theta_0\bigr)'\bigr]
M_n\bigl({\hat{\theta}_{n}^{(J)}}-\theta_0\bigr)\\
&&\qquad\quad{}+\tfrac{1}{2} N^{-1}[E_d\otimes({\hat{\theta}_{n}}-\theta_0)'] T_n({\hat
{\theta}_{n}}-\theta_0)\nonumber\\
&&\qquad\quad{}-N^{-1}\tilde{\Delta}_{n1}\bigl({\hat{\theta}^{(J)}_{n}},
\theta_0\bigr)+N^{-1}\tilde{\Delta}_{n2}({\hat{\theta}_{n}},
\theta_0),
\nonumber
\end{eqnarray}
where $\tilde{\Delta}_{n 1}({\hat{\theta}^{(J)}_{n}}, \theta_0)$ and
$\tilde{\Delta}_{n 2}({\hat{\theta}_{n}}, \theta_0)$ are remainder
terms. Using the same method in the proof of Proposition \ref{prop1}, it can
be shown that $\tilde{\Delta}_{n 1}({\hat{\theta}^{(J)}_{n}},
\theta_0)=O_p\{\|{\hat{\theta}^{(J)}_{n}}- \theta_0\|_2^3\}$ and
$\tilde{\Delta}_{n 2}({\hat{\theta}_{n}},
\theta_0)=O_p\{\|{\hat{\theta}_{n}}-\theta_0\|_2^3\}$.

In order to make $\hat{\theta}_{n}^{(J)}$ have the same asymptotic
distribution as $\hat{\theta}_{n}$, the two quadratic terms on the
right-hand side of (\ref{eqtaylor2}) have to be smaller order of
$\hat{\theta}_{n}^{(J)}-\theta_0$ and $\hat{\theta}_{n} - \theta_0$,
respectively, namely
\[
N^{-1}\bigl[E_d\otimes\bigl({\hat{\theta}_{n}^{(J)}}-\theta_0\bigr)'\bigr]
M_n\bigl({\hat{\theta}_{n}^{(J)}}-\theta_0\bigr)=o_p\bigl\{\bigl\|\hat{\theta
}_{n}^{(J)}-\theta_0\bigr\|_2\bigr\}
\]
or equivalently
%
%e5.2 ###
%e5.2 #&#
\begin{equation}\label{eqcond1}
N^{-1}\bigl[E_d\otimes
\bigl({\hat{\theta}_{n}^{(J)}}-\theta_0\bigr)'\bigr]=o_p(1)
\end{equation}
and
\[
N^{-1}[E_d\otimes({\hat{\theta}_{n}}-\theta_0)']
T_n({\hat{\theta}_{n}}-\theta_0)=o_p\{\|\hat{\theta}_{n}-\theta_0\|_2\}\vadjust{\goodbreak}
\]
or equivalently %$N^{-1}[E_d\otimes(\hat{\theta}_n-\theta_0)^T] =
%o_p(1)$ which is in turn equivalent to
%
%e5.3 ###
%e5.3 #&#
\begin{equation}\label{eqcond2}
n\delta^3\rightarrow\infty
\end{equation}
since ${\hat{\theta}_{n}} - \theta_0 = O_p\{(n\delta)^{-1/2}\}$ and
$N^{-1} = O(\delta^{-1})$.

As ${\hat{\theta}_{n}^{(J)}}-\theta_0 = O_p\{\delta^J +
(n\delta)^{-1/2}\}$, (\ref{eqcond1}) requires that $\delta^{J-1} +
n^{-1/2} \delta^{-3/2} \to0$. Hence, in order to make
$\hat{\theta}_{n}^{(J)}$ have the same asymptotic distribution as~%
$\hat{\theta}_{n}$, it is necessary to have
%
%e5.4 ###
%e5.4 #&#
\begin{equation}\label{eqcond3}
J \geq2\quad\mbox{and}\quad n\delta^3 \to\infty.
\end{equation}
%
%It can be
%readily checked from Proposition 6 and (\ref{eqtaylor2}) that
%under (\ref{eqcond3}),
%$$\sqrt{n} I^{1/2}(\delta) ( {\hat{\theta}_{n}^{(J)}} - \theta_0 )

Now we consider the case of $J=1$. To ensure the remainder terms
$N^{-1}\times\Delta_{n1}({\hat{\theta}_{n}^{(J)}},\theta_0)$ and
$N^{-1}\Delta_{n2}({\hat{\theta}_{n},\theta_0})$ are negligible, by
a similar argument applied above for the case of $J \geq2$, it
is also necessary to assume $n \delta^3 \to\infty$.
%in {\bf under the condition $n\delta^3\rightarrow\infty$}.\footnote{Is
%it too early to talk about
%the need for $n\delta^3\rightarrow\infty$ ? What are the reasons ? \bf
%$n\delta^3\rightarrow\infty$ is necessary to guarantee the quadratic
%term of $\hat{\theta}_n-\theta_0$ is smaller than $\hat{\theta}_n-
%$J=1$?}
From Theorem \ref{theo2}, ${\hat{\theta}_{n}^{(1)}}-\theta_0=O_p\{\delta+
(n\delta)^{-1/2}\}$. To gain insight on the situation, we need to
find out the order of magnitude of the quadratic term in~(\ref{eqtaylor2}), namely the order of magnitude of %which we denote
%as $S_n$ where
%
\[
S_n=N^{-1}\bigl[E_d\otimes\bigl({\hat{\theta}_{n}^{(1)}}-\theta_0\bigr)'\bigr]M_n\bigl({\hat
{\theta}_{n}^{(1)}}-\theta_0\bigr)
-N^{-1}[E_d\otimes({\hat{\theta}_{n}}-\theta_0)']T_n({\hat{\theta
}_{n}}-\theta_0).
\]
With this notation, (\ref{eqtaylor2}) can be written as
%
%e5.5 ###
%e5.5 #&#
\begin{eqnarray}\label{eqtaylor3}
{\hat{\theta}_{n}^{(J)}}-\theta_0&=&N^{-1}U_n+N^{-1}(N_n+F_n)({\hat
{\theta}_{n}}-\theta_0) -\tfrac{1}{2} S_n \nonumber\\[-8pt]\\[-8pt]
&&{}+ o_p\{(n\delta)^{-1/2}\} + O_p(\delta^2).
\nonumber
\end{eqnarray}

Define an operator between two vectors $A$ and $B$,
\[
A*B=[E_d\otimes A']M_nB+[E_d\otimes B'] M_n A.
\]
By repeated substitutions, it can be shown that
\begin{eqnarray*}
S_n&=&%=%&\frac{1}{2} N^{-1}[(N^{-1}U_n)*(N^{-1}U_n)]+\frac{1}{2}
%N^{-1}[(\frac{1}{2} S_n)*(\frac{1}{2} S_n)]-N^{-1}[(N^{-1}U_n)*(
%&+\frac{1}{2} N^{-1}\{[N^{-1}(N_n+F_n)(\hat{\theta}_n-
%&+\frac{1}{2} N^{-1}\{[N^{-1}(N_n-N)(\hat{\theta}_n^{(1)}-
%&+\frac{1}{2} N^{-1}[(N^{-1}\tilde{\Delta}_{n1}-N^{-1}\tilde{
%&+N^{-1}\{(N^{-1}U_n)*[N^{-1}(N_n+F_n)(\hat{\theta}_n-\theta_0)]
%&+N^{-1}[(N^{-1}U_n)*(-N^{-1}\tilde{\Delta}_{n1}+N^{-1}\tilde{
%&-N^{-1}\{[N^{-1}(N_n+F_n)(\hat{\theta}_n-\theta_0)]*[N^{-1}(N_n-N)(
%&-N^{-1}\{[N^{-1}(N_n+F_n)(\hat{\theta}_n-\theta_0)]*(\frac{1}{2} S_n)
%S_n)\}\\
%&+N^{-1}\{[N^{-1}(N_n+F_n)(\hat{\theta}_n-\theta_0)]*(-N^{-1}\tilde{
%&+N^{-1}\{[N^{-1}(N_n-N)(\hat{\theta}_n^{(1)}-\theta_0)]*(-N^{-1}
%&-N^{-1}\{(\frac{1}{2} S_n)*(-N^{-1}\tilde{\Delta}_{n1}+N^{-1}\tilde{
\tfrac{1}{2} N^{-1}[(N^{-1}U_n)*(N^{-1}U_n)]+\tfrac{1}{2} N^{-1}\bigl[\bigl(\tfrac{1}{2}
S_n\bigr)*\bigl(\tfrac{1}{2} S_n\bigr)\bigr]\\
&&{}-N^{-1}\bigl[(N^{-1}U_n)*\bigl(\tfrac{1}{2} S_n\bigr)\bigr]+o_p(\delta).
\end{eqnarray*}

As $U_n= O_p(\delta^2)$ for $J=1$ and $N^{-1} = O(\delta^{-1})$,
it can be deduced from the above equation that $S_n=O_p(\delta)$.
Hence, for $J=1$ if we require $n\delta^3\rightarrow\infty$,
the quadratic term $S_n$ will contribute to the leading order of
$\hat{\theta}_{n}^{(1)}-\theta_0$. If we do not require
$n\delta^3\rightarrow\infty$, then the sum of remainder terms,
$N^{-1}\tilde{\Delta}_{n1}({\hat{\theta}^{(J)}_{n}},
\theta_0)+N^{-1}\tilde{\Delta}_{n2}({\hat{\theta}_{n}}, \theta_0)$
will not be controlled. Hence, if $J=1$, it is very likely that the
asymptotic distribution of $\hat{\theta}_{n}^{(J)}$ will differ from
that of $\hat{\theta}_{n}$ unless $U_n =0$ with probability one. In
the rare case of $U_n = 0$, it is possible for
$\hat{\theta}_{n}^{(1)}$ and~$\hat{\theta}_{n}$ to share the same
limiting distribution.

Therefore, in order to guarantee that $\hat{\theta}_{n}^{(J)}$ has
the same asymptotic distribution as $\hat{\theta}_{n}$ under
$\delta\rightarrow0$, we need to use the AMLE based on at least
two-term expansions, while satisfying
%$J\geq2$ and
$n\delta^3\rightarrow\infty$, which we will assume in the rest of
this section.\vadjust{\goodbreak}%\vspace*{1pt}
%Getting back to $J\geq2$ and $n\delta^3\rightarrow\infty$.

Note that
${\hat{\theta}_{n}^{(J)}}-\theta_0=O_p\{\delta^{J}+(n\delta)^{-1/2}\}$.
Then,
\begin{eqnarray*}
{\hat{\theta}_{n}^{(J)}}-\theta_0
%=&N^{-1}U_n+N^{-1}(N_n+F_n)({\bf\hat{\theta}_{n}}-
%O_p(\delta^{2J}+n^{-1}\delta^{-1})\\
%=&N^{-1}U_n-N^{-1}I({\bf\hat{\theta}_{n}}-\theta_0)+O_p(n^{-1/2}
%=&N^{-1}U_n+(\hat{\theta}_n-\theta_0)+O_p(n^{-1/2}
&=&N^{-1}U_n+({\hat{\theta}_{n}}-\theta_0)\\
&&{}+O_p(n^{-1/2}\delta^{J-3/2})+N^{-1}\cdot O_p(\delta^{2J}+n^{-1}\delta
^{-1}).
\end{eqnarray*}
Furthermore,
\begin{eqnarray*}
&&\sqrt{n}I^{1/2}(\delta)\bigl({\hat{\theta}_{n}^{(J)}}-\theta_0\bigr)\\
&&\qquad=\sqrt{n}I^{-1/2}(\delta)I(\delta) N^{-1}U_n+\sqrt{n}I^{1/2}(\delta)
({\hat{\theta}_{n}}-\theta_0)+O_p(\delta^{J-3/2})\\
&&\qquad\quad{}+\sqrt{n}I^{-1/2}(\delta)I(\delta)N^{-1}\cdot O_p(\delta
^{2J}+n^{-1}\delta^{-1})\\
&&\qquad=\sqrt{n}I^{1/2}(\delta)({\hat{\theta}_{n}}-\theta_0)+O_p(\delta
^{J-3/2}+n^{-1/2}\delta^{-3/2}+n^{1/2}\delta^{J+1/2}).
\end{eqnarray*}
Hence, for any $J\geq2$ such that $n\delta^3\rightarrow\infty$
and $n\delta^{2J+1}\rightarrow0$,
\[
\sqrt{n}I^{1/2}(\delta)\bigl({\hat{\theta}_{n}^{(J)}}-\theta_0\bigr)
\stackrel{d}{\rightarrow}N(0,E_d).
\]
This result shows that, when $\delta$ vanishes to zero, in order to
guarantee the AMLE has the same asymptotic distribution as full MLE,
we need
to pick the approximation order $J\geq2$, while maintaining
$n\delta^3\rightarrow\infty$ and $n\delta^{2J+1}\rightarrow0$. %From
%above description, we
%can get the following theorem about the asymptotic distribution of
%AMLE.

The following theorem summarizes the asymptotic normality under both
asymptotic regimes.
%
%th3 #&#
\begin{tm}\label{theo3}
Under conditions \textup{(A.1)--(A.7)} given in the \hyperref[app]{Appendix},
\[
\sqrt{n}I^{1/2}(\delta)\bigl({\hat{\theta}_{n}^{(J)}}-\theta_0\bigr)
\stackrel{d}{\rightarrow}N(0,E_d)
\]
for: \textup{(i)} $\delta\in(0,\tilde{\Delta}\wedge\dot{\Delta}]$ being
fixed, $n\rightarrow\infty$, $J\rightarrow\infty$ but
$n\delta^{2J+2}\rightarrow0$ or \textup{(ii)}~$J\geq2$ being fixed,
$n\rightarrow\infty$, $\delta\rightarrow0$ but
$n\delta^3\rightarrow\infty$ and $n\delta^{2J+1}\rightarrow0$.
\end{tm}

%s5.3 ###
%s5.3 #&#
\subsection{Asymptotic bias and variance}\label{sec53}

The remainder of this section is devoted to the consideration of the
asymptotic bias and variance of the AMLE under the two asymptotic
regimes. Given our analysis in the early part of this section, our
consideration will be focused on the situations where the asymptotic
normality of the AMLE can be assumed, namely under: (i) $\delta$
being fixed, $J \to\infty$, $n\to\infty$ but $n\delta^{2J+2}\to0$
or (ii) $J \geq2 $ being fixed, $\delta\to0$, $n \delta^3 \to
\infty$ but $n\delta^{2J+1}\to0$.
%which guarantee the AMLE has the same
%asymptotic distribution as MLE.
%By reviewing our analysis made earlier, the asymptotic bias and
%variance for the AMLEs are

In the case of $\delta$ being fixed and $J\rightarrow\infty$, %In 4.1,
%we have considered the asymptotic distribution of AMLE.
%By reviewing our analysis made earlier, provided
%From that part, we know that, under the condition $\delta$ fixed and
%$n\delta^{2J+2}\rightarrow0$,
%the AMLE $\hat{\theta}_n^{(J)}$ has the same asymptotic
%distribution as MLE. Here, we still assume
%$n\delta^{2J+2}\rightarrow0$, then
%$\hat{\theta}_n^{(J)}-\theta_0=O_p(n^{-1/2})$.
%
from (\ref{eqtaylor2}) and provided $n\delta^{2J+2}\rightarrow0$,
we have
\begin{eqnarray*}
{\hat{\theta}_{n}^{(J)}}-\theta_0
&=&N^{-1}U_n+N^{-1}(N_n+F_n)({\hat{\theta}_{n}}-\theta
_0)-N^{-1}(N_n-N)N^{-1}U_n\\
&&{}-N^{-1}(N_n-N)N^{-1}(N_n+F_n)({\hat{\theta}_{n}}-\theta_0)\\
&&{}-\tfrac{1}{2} N^{-1}\{E_d\otimes[N^{-1}U_n+N^{-1}(N_n+F_n)({\hat{\theta
}_{n}}-\theta_0)]'\}\\
&&\hspace*{11.2pt}{}\times M_n[N^{-1}U_n+N^{-1}(N_n+F_n)({\hat{\theta}_{n}}-\theta_0)]\\
&&{}+\tfrac{1}{2} N^{-1}[E_d\otimes({\hat{\theta}_{n}}-\theta_0)']T_n({\hat{\theta
}_{n}}-\theta_0)+O_p(n^{-3/2})\\
&=&N^{-1}U_n+[E_d-N^{-1}(N_n-N)]N^{-1}(N_n+F_n)({\hat{\theta}_{n}}-\theta
_0)\\
&&{}+O_p(n^{-1/2}\delta^{J+1})+O_p(n^{-3/2}).
%=&N^{-1}U_n+[E_d-N^{-1}(N_n-N)]N^{-1}(N_n+F_n)(\hat{\theta}_n-
\end{eqnarray*}
Then, the leading order bias of $\hat{\theta}_{n}^{(J)}$ is
%
%e5.6 ###
%e5.6 #&#
\begin{equation}\label{eqasybias}
B(\theta_0,J,\delta)\,{=}\,N^{-1}U\,{+}\,\mathbb{E}\{[E_d\,{-}\,N^{-1}(N_n\,{-}\,N)]N^{-1}(N_n\,{+}\,F_n)({\hat{\theta}_{n}}\,{-}\,\theta_0)\},\hspace*{-35pt}
\end{equation}
and the leading order variance is
%
%e5.7 ###
%e5.7 #&#
\begin{equation}\label{eqasyvar}
V(\theta_0,J,\delta)=N^{-1}I(\delta)
\operatorname{Var}({\hat{\theta}_{n}})I(\delta) N^{-1}.
\end{equation}

In the case of {$J \geq2$ being fixed, $\delta\rightarrow0$
and $n \delta^3 \to\infty$ but $n \delta^{2J+1} \to0$, it can be
shown by a similar argument to that for the fixed $\delta$ case
above, the asymptotic bias and variance have the same forms as
(\ref{eqasybias}) and (\ref{eqasyvar}), respectively. Both
(\ref{eqasybias}) and (\ref{eqasyvar}) will be used to calibrate
with the simulated bias and variance in the simulation study in
Section \ref{sec7}. For $J=1$ and $\delta\to0$, there are difficulties in
obtaining an expression for the bias {of} the AMLE in general
due to the same dilemma in controlling the reminder terms and the
quadratic term $S_n$ as outlined in Section~\ref{sec52}.

%s6 ###
%s6 #&#
\section{Approximating Fisher information matrix}\label{sec6}

We demonstrate in this section that the approximation of the
transition density provides a way to approximate the Fisher
information matrix. Fisher information matrix~$I(\delta)$ is a key
quantity associated with inference based on the full MLE. It defines
the asymptotic efficiency and {convergence rate}.
%{\bf From Theorem \ref{theo3}, we can find the sense to do the statistical
%inference for the true parameters if we get a reasonable and
%computable substitute to instead of the Fisher information matrix in
%Theorem \ref{theo3}. }
From Proposition~\ref{prop2}, a~natural candidate to approximate
$I(\delta)$ is $-N( \theta_0, J,\delta)$ based on the $J$-term
expansion.
%{\bf However, there are only finite observations in the
%real data set and we do not know the true parameter. Hence, it is
%natural to instead $I(\delta)$ by
%$-N_n(\hat{\theta}_n^{(J)},J,\delta)$.}
To simplify our expedition, our consideration here is
focused under the following diffusion process:
%
%e6.1 ###
%e6.1 #&#
\begin{equation}\label{eq51}
dX_t=\mu(X_t;\eta)\,dt+\sigma(X_t;\xi)\,dB_t,
\end{equation}
where $\eta= (\eta_1, \ldots, \eta_{d_1})'$ and $\xi=(\xi_1, \ldots,
\xi_{d_2})'$ are distinct drift and diffusion parameters,
respectively. The whole parameter $\theta=(\eta',\xi')'$.
Here, we provide an explicit expression $N(\theta_0, 1, \delta)$ based
on the one-term density expansion. Expressions
for higher $J$ values may be made via more extensive derivations.

%we hope to get the first $J$ terms of Taylor expansion of the Fisher
%information by $\delta$, we only need to expand
%$N(\theta_0,J,\delta)$ until $\delta^J$.
Let $\mu_i$, $\mu_{ij}$ and so on denote partial derivatives with
respect to $\eta_i$, $\eta_{i}$ and~$\eta_j$, respectively; and~%
$\sigma_i$ and $\sigma_{x, j}$ and so on denote partial derivatives
with respect to~$\xi_i$, and $x$ and $\xi_j$, respectively. By the
one-term ($J=1$) transition density approximation, derivations given
in \citet{r14} show that
\[
\mathbb{E}\biggl(\frac{\partial^2\log
f^{(1)}}{\partial\eta_i\,\partial\eta_j}\biggr)=:\delta\cdot
N_{11}^{(1)} +O(\delta^2),\qquad \mathbb{E}\biggl(\frac{\partial^2\log
f^{(1)}}{\partial\eta_i\,\partial\xi_j}\biggr)=:\delta\cdot
N_{12}^{(1)} +O(\delta^2)
\]
and
\[
\mathbb{E}\biggl(\frac{\partial^2\log
f^{(1)}}{\partial\xi_i\,\partial\xi_j}\biggr)
=:-2\mathbb{E}(\sigma^{-2}\sigma_i\sigma_j) + \delta\cdot
N^{(1)}_{22} + O(\delta^2),
\]
where
\begin{eqnarray*}
N_{11}^{(1)} &=& \mathbb{E}\bigl\{-\sigma^{-2}\mu_i\mu_j-\mu\sigma^{-2}\mu
_{ij}+\sigma^{-1}\mu_{ij}\sigma_x-\tfrac{1}{2}\mu_{xij}\bigr\},\\
N_{12}^{(1)} &=& \mathbb{E}\{2\mu\sigma^{-3}\mu_i\sigma_j-\sigma^{-2}\mu
_i\sigma_x\sigma_j+\sigma^{-1}\mu_i\sigma_{xj}\},\\
N^{(1)}_{22} &=& \mathbb{E}\bigl\{-6\mu^2\sigma^{-4}\sigma_i\sigma_j+16\mu
\sigma^{-3}\sigma_x\sigma_i\sigma_j+2\mu^2\sigma^{-3}\sigma_{ij}-3\sigma
^{-2}\mu_x\sigma_i\sigma_j\\
&&\hphantom{\mathbb{E}\bigl\{}{}
-\tfrac{19}{2}\sigma^{-2}\sigma_x^2\sigma_i\sigma_j-\tfrac{9}{2}\mu
\sigma^{-2}\sigma_x\sigma_{ij}-5 \mu\sigma^{-2}\sigma_{xi}\sigma_j -
5\mu\sigma^{-2}\sigma_{xj}\sigma_i\\
&&\hphantom{\mathbb{E}\bigl\{}{}
+ \sigma^{-1}\mu_x\sigma_{ij}+4\sigma^{-1}\sigma_{xx}\sigma_i\sigma
_j+\tfrac{11}{2}\sigma^{-1}\sigma_x\sigma_{xi}\sigma_j+\tfrac{11}{
2}\sigma^{-1}\sigma_x\sigma_{xj}\sigma_i\\
&&\hphantom{\mathbb{E}\bigl\{}{}
+ \tfrac{3}{2}\sigma^{-1}\sigma_x^2\sigma_{ij}+\tfrac{5}{2}\mu\sigma
^{-1}\sigma_{xij}-\tfrac{3}{4}\sigma_{xx}\sigma_{ij}-\tfrac{5}{
2}\sigma_{xi}\sigma_{xj}-\tfrac{3}{2}\sigma_x\sigma_{xij}\\
&&\hphantom{\mathbb{E}\bigl\{}\hspace*{153pt}{}
-\sigma_{xxi}\sigma_j-\sigma_{xxj}\sigma_i+\tfrac{3}{4}\sigma\sigma
_{xxij}\bigr\}.
\end{eqnarray*}
Thus
%
%e6.2 ###
%e6.2 #&#
\begin{equation}\label{eqfisherapprox}
N(\theta_0,1,\delta) = \pmatrix{
\delta\cdot N_{11}^{(1)} & \delta\cdot N_{12}^{(1)}
\vspace*{2pt}\cr
\delta\cdot{N_{12}^{(1)}}' & -2\cdot\mathbb{E}(\sigma^{-2}\sigma
_i\sigma_j) + \delta\cdot N^{(1)}_{22}}
 + O(\delta^2).
\end{equation}

We learn from Proposition \ref{prop2} that $-N(\theta_0,1,\delta)$ provides a
leading order approximation to $I(\delta)$ with a reminder term at
the order of $\delta^2$. Equation~(\ref{eqfisherapprox}) confirms
that as $\delta\to0$, given the asymptotic normality of the full
MLE $\hat{\theta}_{n}$ as conveyed by Proposition \ref{prop6}, that the
convergence rate of the full MLE for the drift parameters $\eta$ is
$(n\delta)^{-1/2}$ whereas that for the diffusion parameters~$\xi$
is $n^{-1/2}$, faster than the drift parameter estimator. Our study
confirms the results of \citet{r23}, \citet{r32} and \citet{r36}. %The result in
%(\ref{eqfisherapprox}) provides extra reasoning the differential
%convergence rates
%between the drift and diffusion parameters for much general
%diffusion processes.

In the rest of the section, we will derive the Fisher information
matrix approximation for two specific diffusion processes. Both are
widely employed in {modeling} of the interest rate dynamics.

%s6.1 ###
%s6.1 #&#
\subsection{Vasicek model}\label{sec61}

Consider the Vasicek (\citeyear{r37}) model,
%
%e6.3 ###
%e6.3 #&#
\begin{equation}
dX_t=\kappa(\alpha-X_t)\,dt+\sigma \,dB_t,
\end{equation}
which is also the Ornstein--Uhlenbeck process. The conditional
distribution of $X_t$ given $X_{t-1}$ is
\[
X_t|X_{t-1}\sim
N\bigl\{X_{t-1}e^{-\kappa\delta}+\alpha(1-e^{-\kappa\delta}),\tfrac{1}{2}
\sigma^2\kappa^{-1}(1-e^{-2\kappa\delta})\bigr\},%\label{eqvtransitional}
\]
and the stationary distribution of $\{X_t\}$ is $
N(\alpha,\frac{\sigma^2}{2\kappa}).\label{eqvasicekstationarydistribution}
$
It yields that the information matrix of
$\theta=(\kappa,\alpha,\sigma)'$ is $I(\delta)=(
I_{ij})_{3 \times3}$
%I(\delta, \theta)=(
% \begin{array}{ccc}
% I_{11}&I_{12}& I_{13}\\
% I_{21} & I_{22} & I_{23} \\
% I_{31} & I_{32} & I_{33} \\
% \end{array}
%) \label{eqfisherinformation}
where
\begin{eqnarray*}
I_{11}&=&\frac{1}{2\kappa^2}+\frac{\delta[\kappa\delta+\kappa\delta
e^{2\kappa\delta}-2e^{2\kappa\delta}+2]}{\kappa(e^{2\kappa\delta
}-1)^2}=\frac{\delta}{2\kappa}+O(\delta^2),\qquad I_{12}=I_{21}=0,
\\
I_{13}&=&I_{31}=\frac{(1+2\kappa\delta)-e^{2\kappa\delta}}{\sigma\kappa
(e^{2\kappa\delta}-1)}=-\frac{\delta}{\sigma}+O(\delta^2),
\\
I_{22}&=&\frac{2\kappa(e^{\kappa\delta}-1)^2}{\sigma^2(e^{2\kappa\delta
}-1)}=\frac{\kappa^2\delta}{\sigma^2}+O(\delta^2),\qquad
I_{23}=I_{32}=0\quad\mbox{and}\quad I_{33}=\frac{2}{\sigma^2}.
\end{eqnarray*}
These mean that
%
%e6.4 ###
%e6.4 #&#
\begin{equation} \label{eqfishervas}
I(\delta) = \pmatrix{
\delta\cdot(2\kappa)^{-1} & 0 & -\delta\cdot\sigma^{-1} \cr
0 & \delta\cdot\kappa^2\sigma^{-2} & 0 \cr
-\delta\cdot\sigma^{-1} & 0 & 2\sigma^{-2}}
+O(\delta^2).
\end{equation}
Hence $I(0) =
\lim_{\delta\to0} I(\delta)$ is singular, an issue we have raised
earlier, which makes us assume that $\delta I^{-1}(\delta)$'s largest
eigenvalue is bounded in condition~(A.5).

Using the approximation formula in (\ref{eqfisherapprox}), we
have
\[
N(\theta, 1, \delta) = \pmatrix{
-\delta\cdot(2\kappa)^{-1} & 0 & \delta\cdot\sigma^{-1} \cr
0 & -\delta\cdot\kappa^2\sigma^{-2} & 0 \cr
\delta\cdot\sigma^{-1} & 0 & -2\sigma^{-2} }+O(\delta^2).
\]
This means the leading order term of $-N(\theta, 1, \delta)$ is
identical with that of the true Fisher information matrix in
(\ref{eqfishervas}).

%s6.2 ###
%s6.2 #&#
\subsection{Cox--Ingersoll--Ross model}\label{sec62}
Consider the Cox--Ingersoll--Ross (CIR) model
[Cox, Ingersoll and Ross (\citeyear{r15})],
%
%e6.5 ###
%e6.5 #&#
\begin{equation}
dX_t=\kappa(\alpha-X_t)\,dt+\sigma\sqrt{X_t}\,dB_t,
\end{equation}
which is also Feller's (\citeyear{r21}) square root {process}.

Let $\theta=(\kappa,\alpha,\sigma)'$ and
$c=4\kappa\sigma^{-2}(1-e^{-\kappa\delta})^{-1}$. The conditional
distribution of $cX_t$ given $X_{t-1}$ is
\[
cX_t|X_{t-1}\sim\chi^2_\nu(\lambda),
\]
where the distribution is a noncentral $\chi^2$ distribution with
degree of freedom $\nu=4\kappa\alpha\sigma^{-2}$ and noncentral
parameter $\lambda=cX_{t-1}e^{-\kappa\delta}$. The transition
density of $X_{t}$ given $X_{t-1}$ is
\[
f(X_t|X_{t-1},\delta;\theta)=\frac{c}{2}e^{-u-v}\biggl(\frac
{v}{u}\biggr)^{q/2}I_q\bigl(2\sqrt{uv}\bigr),
\]
where $u=cX_{t-1}e^{-\kappa\delta}/2$, $v=cX_t/2$,
$q=2\kappa\alpha/\sigma^2-1\geq0$, and $I_q$ is the modified
Bessel function of the first kind of order $q$. If
$2\kappa\alpha>\sigma^2$, then the stationary distribution of
$\{X_t\}$ is
$\Gamma(\frac{2\kappa\alpha}{\sigma^2},\frac{\sigma^2}{2\kappa}).\label
{eqstationary}
$
%The log {\bf transition} density function is
%v-\log u)+\log I_q(2\sqrt{uv})-\log2.

Although the second partial derivations of the log transition
density function can be derived after some labor that is involved with
differentiating the modified Bessel function of {the} first kind,
acquiring an expression for the Fisher information matrix is a
rather hard task, largely due to the difficulty in deriving the
expectations.
%We can use the second derivatives of log transitional density function
%to get the Fisher information matrix of CIR model. But it's very hard
%to get the exact expression of each element in the matrix, as the
%expressions of second derivatives contain $\partial
%I_l(x)/\partial l$ and $\partial^2I_l(x)/\partial l^2$, which are too
%hard to be simplified.
In contrast, using the approximation formula~(\ref{eqfisherapprox}), we can obtain the approximation for the
opposite Fisher information matrix,
\[
N(\theta_0, 1, \delta) = \pmatrix{
N_{11} & N_{12} & N_{13} \cr
N_{21} & N_{22} & N_{23} \cr
N_{31} & N_{32} & N_{33} }
+O(\delta^2),
\]
where
%N_{11}=\delta\cdot\sigma^{-2}\cdot\mathbb{E}\{X_t^{-1}(\alpha-X_t)^2
%N_{13}=N_{31}=-\delta\cdot2\kappa\sigma^{-3}\cdot\mathbb{E}\{X_t^{-1}(
%N_{23}=N_{32}=-\delta\cdot2\kappa^2\sigma^{-3}\cdot\mathbb{E}
%N_{33} = 2\sigma^{-2}
%-\delta\cdot3\kappa\sigma^{-2}+\delta\cdot\mathbb{E}\{6\kappa^2
%More explicit form of the approximation may be obtained by
%cultivating the marginal distribution of $X_t$. Under
%(\ref{eqstationary}), we can get
%Then,
%
\begin{eqnarray*}
N_{11}&=&\delta\cdot\frac{\alpha^2\sigma^2-2\kappa\alpha^2+\alpha\sigma
^2}{\sigma^4-2\kappa\alpha\sigma^2},\\
N_{12}&=&N_{21}=\delta\cdot\frac
{4\kappa\alpha\sigma^2-\sigma^4-8\kappa^2\alpha+4\kappa\sigma^2}{2\sigma
^4-4\kappa\alpha\sigma^2},\\
N_{13}&=&N_{31}=-\delta\cdot\frac{2\kappa\alpha^2\sigma^2-4\kappa^2\alpha
^2+2\kappa\alpha\sigma^2}{\sigma^5-2\kappa\alpha\sigma^3},\\
N_{22}&=&\delta\cdot\frac{\kappa^2}{\sigma^2-2\kappa\alpha},\\
N_{23}&=&-\delta\cdot\frac{2\kappa^2\alpha\sigma^2-4\kappa^3\alpha+2\kappa
^2\sigma^2}{\sigma^5-2\kappa\alpha\sigma^3}
\end{eqnarray*}
and
\begin{eqnarray*}
N_{33}&=&\frac{-2}{\sigma^2}
+\delta
\cdot\fracb{24\kappa^2\alpha^2\sigma^2-48\kappa^3\alpha^2+48\kappa
^2\alpha\sigma^2\\
&&\hspace*{44pt}{}-24\kappa\alpha\sigma^4+36\kappa\sigma^4+4\sigma
^5+9\sigma^6}{4\sigma^6-8\kappa\alpha\sigma^4}.
\end{eqnarray*}

Using $-N(\theta_0, 1, \delta)$, we can get the approximation of the
Fisher information matrix. This approximation may be used in
carrying out statistical inference on the CIR processes.

%s6.3 ###
%s6.3 #&#
\subsection{Observed Fisher information}\label{sec63}

The major application for the asymptotic normality of both the full
and approximate MLEs is for statistical inference of $\theta$, which
include confidence regions and testing hypotheses for $\theta$. For
such purposes, the Fisher information $I(\delta)$ needs to be
estimated. A natural candidate would be $-N_n(\hat{\theta}_n^{(J)},
J, \delta)$. Although\vspace*{1pt} it converges to $I(\delta)$ at the rate of
$O_p\{ (n\delta)^{-1/2} + \delta^J\}$ or $O_p\{ (n\delta)^{-1/2} +
\delta^{J+1}\}$, depending on whether $\delta$ is fixed or diminishing,
$-N_n(\hat{\theta}_n^{(J)}, J, \delta)$ may not be nonnegative
definite, which can hinder the acquisition of
$\{-N_n(\hat{\theta}_n^{(J)}, J, \delta)\}^{1/2}$. To get around
this issue, by noticing that $I(\delta)$ is the variance of the
likelihood score, we consider
\[
\tilde{I}_n(\theta,J,\delta)=\frac{1}{n}\sum_{t=1}^n\bigl[\nabla_\theta
\log
f^{(J)}(X_t|X_{t-1},\delta;\theta)\bigr]\bigl[\nabla_\theta\log
f^{(J)}(X_t|X_{t-1},\delta;\theta)\bigr]'
\]
as an estimator of $I(\delta)$. The following theorem shows this by
replacing~$I(\delta)$ with
$\tilde{I}_n(\hat{\theta}_n^{(J)},J,\delta)$ in Theorem \ref{theo3}.
%
%th4 #&#
\begin{tm}\label{theo4}
Under conditions \textup{(A.1)--(A.7)} given in the \hyperref[app]{Appendix},
\[
\sqrt{n}\tilde{I}_n^{1/2}\bigl(\hat{\theta}_n^{(J)},J,\delta\bigr)\bigl(\hat{\theta
}_n^{(J)}-\theta_0\bigr)
\stackrel{d}{\rightarrow}N(0,E_d)
\]
for: \textup{(i)} $\delta\in(0,\tilde{\Delta}\wedge\dot{\Delta}]$ being
fixed, $n\rightarrow\infty$, $J\rightarrow\infty$ but
$n\delta^{2J+2}\rightarrow0$ or \textup{(ii)}~\mbox{$J\geq2$} being fixed,
$n\rightarrow\infty$, $\delta\rightarrow0$ but
$n\delta^3\rightarrow\infty$ and $n\delta^{2J+1}\rightarrow0$.
\end{tm}

Confidence regions and testing hypotheses can be readily carried out by
utilizing the above results.

%s7 ###
%s7 #&#
\section{Simulation}\label{sec7}

We report results from simulation studies which are designed to
confirm the theoretical findings on the AMLE as reported in the
earlier sections. To allow verification with the full MLE, we
considered the Vasicek and CIR diffusion models reported in the
previous section as both models permit the full MLE. The two
asymptotic regimes were experimented: the fixed $\delta$ and the
diminishing $\delta$
%There are two parts in this section. Firstly, we consider the
%simulation for Vasicek model and CIR model under fixed $\delta$.
%Secondly, we consider the simulation for AMLE with approximation
%order 1 in Vasicek model under $\delta$ vanish to zero, but
with $n\delta^3\rightarrow\infty$.

The first part of the simulation is about the case in which $\delta$ is
fixed. The parameters used in the simulated Vasicek and CIR models
were $\theta=(\kappa,\alpha,\sigma)'=(0.858,0.0891,0.0468)'$ and
$\theta=(\kappa,\alpha,\sigma)'=(0.892,0.09,0.1817)'$, respectively.
The sampling interval $\delta$ was $1/12$ and $1/4$, and the order
of the density approximation $J$ was 1 and 2, respectively. For each
$\delta$ and $J$, the sample size $n$ was set at $500$, 1,000 and
2,000, respectively. In addition to bias and standard deviation, we
consider
\[
\mathrm{RMSD}(n,J,\delta)=\sqrt{\mathbb{E}\bigl\|{\hat{\theta}_{n}^{(J)}-\hat
{\theta}_{n}}\bigr\|^2_2},%\vadjust{\goodbreak}
\]
the\vspace*{1pt} square root of the expected square of modulated deviations
between $\hat{\theta}_{n}^{(J)}$ and~$\hat{\theta}_{n}$, as an overall
performance measure. %compare the performance of AMLE with different
%approximation order. %At the same
%time, in order to find the approximate bias and approximate standard
%variance for $J=1$ and $J=2$, we need to know
%$N^{-1}(\theta,J,\delta)U(\theta,J,\delta)$ and
%$N^{-1}(\theta,J,\delta)\mathbf{I}(\delta)$. We simulate 20000
%samples to get the two.

%
%t2 ###
%t2 #&#
\begin{table}%[b]
\caption{Simulated average bias (Bias) and
standard deviations (SD) of the full MLE and two AMLEs
with $J=1$ and $2$ for Vasicek
model ($\kappa=0.858, \alpha=0.0891, \sigma=0.0468$); A.Bias and
A.SD are asymptotic bias and SD based on formulas (\protect\ref{eqasybias})
and (\protect\ref{eqasyvar}); RMSD is the root of mean square deviation
between $\hat{\theta}_{n}$ and $\hat{\theta}_{n}^{(J)}$}
\label{table2}
\begin{tabular*}{\tablewidth}{@{\extracolsep{\fill}}lccd{2.5}d{2.5}d{2.5}d{1.5}d{2.5}d{1.5}@{}}
\hline
& &&\multicolumn{3}{c}{$\bolds{\delta=1/12}$}
& \multicolumn{3}{c@{}}{$\bolds{\delta=1/4}$}\\[-4pt]
& &&\multicolumn{3}{c}{\hrulefill}&\multicolumn{3}{c@{}}{\hrulefill}\\
$\bolds{n}$
& \multicolumn{2}{c}{\textbf{Statistics}} &
\multicolumn{1}{c}{\textbf{MLE}} & \multicolumn{1}{c}{$\bolds{J=1}$}
& \multicolumn{1}{c}{$\bolds{J=2}$} & \multicolumn{1}{c}{\textbf{MLE}}
& \multicolumn{1}{c}{$\bolds{J=1}$} & \multicolumn{1}{c@{}}{$\bolds{J=2}$}\\
\hline
500& Bias &$\kappa$&0.0992 & 0.0896 &0.0992 &0.0380 &0.0127 & 0.0396 \\
&&$\alpha$&0.0002 &0.0002 &0.0002 &4.09\mbox{e--}5 &5.63\mbox{e--}5 &4.17\mbox{e--}5 \\
& &$\sigma$&4.39\mbox{e--}5 &4.14\mbox{e--}5 &4.39\mbox{e--}5 &9.12\mbox{e--}5 &7.13\mbox{e--}5 &9.43\mbox{e--}5 \\
[3pt]
&A.Bias &$\kappa$& &0.0908 &0.1016 & & 0.0174 &0.0376 \\
&&$\alpha$& &0.0003 &0.0002 & &0.0002 &0.0001 \\
& &$\sigma$& & 4.55\mbox{e--}5 &4.55\mbox{e--}5 & &0.0001 &0.0001 \\
[3pt]
&SD &$\kappa$&0.2307 &0.2255 &0.2309 &0.1366 &0.1290 &0.1386 \\
& &$\alpha$&0.0085 &0.0085 &0.0085 &0.0050 &0.0050 &0.0050 \\
& &$\sigma$&0.0016 &0.0016 &0.0016 &0.0016 &0.0016 &0.0016 \\
[3pt]
&A.SD &$\kappa$& &0.2251 &0.2366 & &0.1215 &0.1403 \\
&&$\alpha$& &0.0084 &0.0085 & &0.0047 &0.0050 \\
& &$\sigma$& &0.0016 &0.0016 & &0.0016 &0.0016 \\
[3pt]
&RMSD &$\kappa$& &0.0173 &0.0062 & &0.0332 &0.0316 \\
&&$\alpha$& &0.0002 &1.28\mbox{e--}5 & &0.0005 &0.0002 \\
& &$\sigma$& &1.36\mbox{e--}5 &1.05\mbox{e--}5 & &0.0001 &0.0001 \\
[6pt]
1,000& Bias&$\kappa$&0.0518 & 0.0419 &0.0520 &0.0170 &-0.0095 &0.0186 \\
&&$\alpha$&-0.0002 &-0.0002 &-0.0002 &1.83\mbox{e--}5 &2.81\mbox{e--}5 &1.58\mbox{e--}5 \\
& &$\sigma$&7.05\mbox{e--}5 &6.68\mbox{e--}5 &7.06\mbox{e--}5 &3.66\mbox{e--}5 &6.83\mbox{e--}6 &3.96\mbox{e--}5 \\
[3pt]
& A.Bias&$\kappa$& & 0.0446 &0.0529 & &-0.0097 &0.0161 \\
&&$\alpha$& &-0.0001 &-0.0002 & &1.69\mbox{e--}5 &1.45\mbox{e--}5 \\
& &$\sigma$& &0.0001 &0.0001 & &3.29\mbox{e--}5 &4.55\mbox{e--}5 \\
[3pt]
&SD &$\kappa$&0.1624 &0.1586 &0.1625 &0.0957 &0.0905 &0.0966 \\
& &$\alpha$&0.0058 &0.0058 &0.0058 &0.0034 &0.0034 &0.0034 \\
& &$\sigma$&0.0011 &0.0011 &0.0011 &0.0012 &0.0012 &0.0012 \\
[3pt]
&A.SD &$\kappa$& &0.1585 &0.1666 & &0.0849 &0.0982 \\
&&$\alpha$& &0.0057 &0.0058 & &0.0032 &0.0034 \\
& &$\sigma$& &0.0011 &0.0011 & &0.0012 &0.0012 \\
[3pt]
&RMSD &$\kappa$& &0.0100 &0.0008 & &0.0316 &0.0063 \\
&&$\alpha$& &0.0001 &9.14\mbox{e--}6 & &0.0004 &0.0001 \\
& &$\sigma$& &7.39\mbox{e--}6 &7.80\mbox{e--}7 & &0.0001 &1.59\mbox{e--}5 \\
\hline
\end{tabular*}
\end{table}
\setcounter{table}{1}
\begin{table}
\caption{(Continued)}
\begin{tabular*}{\tablewidth}{@{\extracolsep{\fill}}lccd{2.5}
d{2.5}d{2.5}d{2.5}d{2.5}d{2.5}@{}}
\hline
& &&\multicolumn{3}{c}{$\bolds{\delta=1/12}$}
& \multicolumn{3}{c@{}}{$\bolds{\delta=1/4}$}\\[-4pt]
& &&\multicolumn{3}{c}{\hrulefill}&\multicolumn{3}{c@{}}{\hrulefill}\\
$\bolds{n}$
& \multicolumn{2}{c}{\textbf{Statistics}} &
\multicolumn{1}{c}{\textbf{MLE}} & \multicolumn{1}{c}{$\bolds{J=1}$}
& \multicolumn{1}{c}{$\bolds{J=2}$} & \multicolumn{1}{c}{\textbf{MLE}}
& \multicolumn{1}{c}{$\bolds{J=1}$} & \multicolumn{1}{c@{}}{$\bolds{J=2}$}\\
\hline
2,000&Bias &$\kappa$&0.0245 &0.0149 &0.0246 &0.0084 &-0.0191 &0.0100 \\
&&$\alpha$&-3.97\mbox{e--}5 &-3.34\mbox{e--}5 &-4.01\mbox{e--}5 &-5.72\mbox{e--}5 &-4.90\mbox{e--}5
&-5.80\mbox{e--}5 \\
& &$\sigma$&2.69\mbox{e--}5 &2.30\mbox{e--}5 &2.70\mbox{e--}5 &4.00\mbox{e--}5 &9.21\mbox{e--}6 &4.34\mbox{e--}5 \\
[3pt]
&A.Bias &$\kappa$& &0.0179 &0.0249 & &-0.0085 & 0.0071 \\
&&$\alpha$& &-2.63\mbox{e--}5 &-2.98\mbox{e--}5 & &0.0001 &-0.0001 \\
& &$\sigma$& &4.55\mbox{e--}5 &4.55\mbox{e--}5 & &4.55\mbox{e--}5 & 4.55\mbox{e--}5 \\
[3pt]
&SD &$\kappa$&0.1114 &0.1091 &0.1115 &0.0647 &0.0611 &0.0652 \\
& &$\alpha$&0.0042 &0.0041 &0.0042 &0.0024 &0.0024 &0.0024 \\
& &$\sigma$&0.0008 &0.0008 &0.0008 &0.0008 &0.0008 &0.0008 \\
[3pt]
&A.SD &$\kappa$& &0.1088 &0.1143 & &0.0576 &0.0665 \\
&&$\alpha$& &0.0041 &0.0042 & &0.0023 &0.0024 \\
& &$\sigma$& &0.0008 &0.0008 & &0.0008 &0.0008 \\
[3pt]
&RMSD &$\kappa$& &0.0100 &0.0006 & &0.0300 &0.0042 \\
&&$\alpha$& &0.0001 &7.37\mbox{e--}6 & &0.0003 &4.68\mbox{e--}5 \\
& &$\sigma$& &6.27\mbox{e--}6 &7.80\mbox{e--}7 & &0.0001 &1.02\mbox{e--}5 \\
\hline
\end{tabular*}
\end{table}

Tables \ref{table2} and \ref{table3} summarize the simulation for the fixed $\delta$ case.
They report the average bias and standard deviation (SD) for the
full MLE and AMLEs with $J=1$ and $J=2$, as well as the RMSD between
the AMLEs and the full MLE, for both the Vasicek and the CIR models.
To give the simulation results more perspective and to confirm the
derived {approximate} bias and variance formulas in Section \ref{sec5}, we
also computed the asymptotic bias and standard deviation based on
formulas (\ref{eqasybias}) and (\ref{eqasyvar}). We observe
from Tables \ref{table2} and \ref{table3} that at each $\delta$ ($1/12$ and $1/4$)
experimented,
the bias and the standard deviation of all {the} estimators for {the}
three parameters became smaller as $n$ increased. These confirmed the
consistency of the estimators. The tables also showed that there was a
good agreement among the three estimators in terms of the performance
measures. It appeared that the bias and the variance of the AMLE
with $J=1$ and $J=2$ were quite comparable to each other. However, by
comparing RMSD, it was clear that in most of the cases (except for
$n=500$ of CIR model), the RMSD for $J=2$ was smaller than $J=1$,
signaling the AMLE with $J=2$ was closer to the full MLE than that of
the AMLE with $J=1$. This indicates that the AMLEs with $J=2$ were
indeed closer to those with $J=1$, as confirmed by our early analysis.
%% and the studies in A\"{\i}t-Sahalia (1999).
The asymptotic bias and standard deviation predicted for the AMLE with
$J=1$ and $2$ offer more insights, and show good agreement between the
simulated results and the predicted values by the theory, which is very
assuring.
We also {observe} that for $\delta=1/4$, the AMLE with $J=2$
performs better than AMLE with $J=1$, which somehow reflects Table \ref{table1}
which shows that $J=2$ is preferred to $J=1$ at this frequency.
When $\delta$ was fixed at $1/12$, we see the performance between
$J=1$ and $J=2$ was largely similar.

%t3 ###
%t3 #&#
\begin{table}
\caption{Simulated average bias (Bias) and standard deviations (SD)
of the full MLE and two AMLEs with $J=1$ and $2$ for CIR model
($\kappa=0.892, \alpha=0.09, \sigma=0.1817$); A.Bias and A.SD are
asymptotic bias and SD based on formulas (\protect\ref{eqasybias}) and
(\protect\ref{eqasyvar}); RMSD is the root of mean square deviation between
$\hat{\theta}_{n}$ and $\hat{\theta}_{n}^{(J)}$}
\label{table3}
\begin{tabular*}{\tablewidth}{@{\extracolsep{\fill}}lccd{2.5}d{1.5}d{2.5}d{2.5}cd{2.4}@{}}
\hline
& &&\multicolumn{3}{c}{$\bolds{\delta=1/12}$}&\multicolumn{3}{c@{}}{$\bolds{\delta
=1/4}$}\\[-4pt]
& &&\multicolumn{3}{c}{\hrulefill}&\multicolumn{3}{c@{}}{\hrulefill}\\
$ \bolds{n} $& \multicolumn
{2}{c}{\textbf{Statistics}}&\multicolumn{1}{c}{\textbf{MLE}}
&\multicolumn{1}{c}{$\bolds{J=1}$}&\multicolumn{1}{c}{$\bolds{J=2}$}
&\multicolumn{1}{c}{\textbf{MLE}}&\multicolumn{1}{c}{$\bolds{J=1}$}
&\multicolumn{1}{c@{}}{$\bolds{J=2}$}\\
\hline
500& Bias&$\kappa$&0.0980 &0.0910 &0.0978 &0.0371 &0.0234 & 0.0388 \\
&&$\alpha$&0.0001 &0.0004 &0.0001 &-6.38\mbox{e--5} &0.0008 &-0.0001 \\
& &$\sigma$&0.0003 &0.0003 &0.0003 &0.0004 &0.0005 &0.0003 \\
[3pt]
&A.Bias &$\kappa$& &0.0818 &0.0984 & &0.0207 & 0.0513 \\
&&$\alpha$& &0.0005 &0.0001 & &0.0008 &-0.0001 \\
& &$\sigma$& &0.0003 &0.0003 & &0.0004 &0.0002 \\
[3pt]
&SD &$\kappa$&0.2389 &0.2340 &0.2405 &0.1437 &0.1338 &0.2256 \\
& &$\alpha$&0.0093 &0.0093 &0.0093 &0.0055 &0.0054 &0.0055 \\
& &$\sigma$&0.0060 &0.0060 &0.0060 &0.0065 &0.0065 &0.0069 \\
[3pt]
&A.SD &$\kappa$& &0.2169 &0.2389 & &0.1159 &0.1938 \\
&&$\alpha$& &0.0091 &0.0093 & &0.0064 &0.0055 \\
& &$\sigma$& &0.0060 &0.0060 & &0.0067 &0.0065 \\
[3pt]
&RMSD &$\kappa$& &0.0200 &0.0224 & &0.0447 &0.1622 \\
&&$\alpha$& &0.0009 &0.0004 & &0.0018 &0.0004 \\
& &$\sigma$& &0.0004 &0.0004 & &0.0017 &0.0021 \\
[6pt]
1,000& Bias&$\kappa$&0.0521 &0.0435 &0.0521 &0.0218 &0.0070 &0.0186 \\
&&$\alpha$&-1.54\mbox{e--}5 &0.0002 &-2.22\mbox{e--}5 &-0.0002 &0.0007 &-0.0003 \\
& &$\sigma$&3.86\mbox{e--}5 &4.35\mbox{e--}5 &3.81\mbox{e--}5 &0.0003 &0.0006 &0.0003 \\
[3pt]
&A.Bias &$\kappa$& &0.0411 &0.0525 & &0.0095 &0.0262 \\
&&$\alpha$& &0.0004 &-3.43\mbox{e--}5 & &0.0007 &-0.0003 \\
& &$\sigma$& &3.17\mbox{e--}5 &2.69\mbox{e--}5 & &0.0003 &0.0001 \\
[3pt]
&SD &$\kappa$&0.1596 &0.1558 &0.1603 &0.0968 &0.0861 &0.0980 \\
& &$\alpha$&0.0067 &0.0067 &0.0067 &0.0039 &0.0037 &0.0039 \\
& &$\sigma$&0.0043 &0.0043 &0.0043 &0.0045 &0.0045 &0.0045 \\
[3pt]
&A.SD &$\kappa$& &0.1452 &0.1596 & &0.0823 &0.0969 \\
&&$\alpha$& &0.0066 &0.0067 & &0.0044 &0.0039 \\
& &$\sigma$& &0.0040 &0.0043 & &0.0047 &0.0045 \\
[3pt]
&RMSD &$\kappa$& &0.0173 &0.0141 & &0.0447 &0.0200 \\
&&$\alpha$& &0.0003 &2.66\mbox{e--}5 & &0.0020 &0.0001 \\
& &$\sigma$& &0.0002 &3.91\mbox{e--}5 & &0.0021 &0.0002 \\
\hline
\end{tabular*}
\vspace*{-5pt}
\end{table}

\setcounter{table}{2}
\begin{table}
\caption{(Continued)}
\begin{tabular*}{\tablewidth}{@{\extracolsep{\fill}}lccd{2.4}cd{2.4}d{2.5}d{2.4}d{2.5}@{}}
\hline
& &&\multicolumn{3}{c}{$\bolds{\delta=1/12}$}&\multicolumn{3}{c@{}}{$\bolds{\delta
=1/4}$}\\[-4pt]
& &&\multicolumn{3}{c}{\hrulefill}&\multicolumn{3}{c@{}}{\hrulefill}\\
$ \bolds{n} $& \multicolumn
{2}{c}{\textbf{Statistics}}&\multicolumn{1}{c}{\textbf{MLE}}
&\multicolumn{1}{c}{$\bolds{J=1}$}&\multicolumn{1}{c}{$\bolds{J=2}$}
&\multicolumn{1}{c}{\textbf{MLE}}&\multicolumn{1}{c}{$\bolds{J=1}$}
&\multicolumn{1}{c@{}}{$\bolds{J=2}$}\\
\hline
2,000&Bias &$\kappa$&0.0295 &0.0199 &0.0294 &0.0103 &-0.0057 &0.0069 \\
&&$\alpha$&-0.0002 &0.0001 &-0.0002 &-3.06\mbox{e--}5 &0.0010 &-9.87\mbox{e--}5 \\
& &$\sigma$&0.0002 &0.0002 &0.0002 &3.05\mbox{e--}5 &0.0006 &1.33\mbox{e--}5 \\
[3pt]
&A.Bias &$\kappa$& &0.0213 &0.0299 & &-0.0011 &0.0147 \\
&&$\alpha$& &0.0002 &-0.0002 & &0.0006 &-0.0001 \\
& &$\sigma$& &0.0002 &0.0002 & &0.0005 &1.06\mbox{e--}5 \\
[3pt]
&SD &$\kappa$&0.1082 &0.1053 &0.1088 &0.0696 &0.0607 &0.0698 \\
& &$\alpha$&0.0048 &0.0048 &0.0048 &0.0028 &0.0027 &0.0028 \\
& &$\sigma$&0.0030 &0.0031 &0.0030 &0.0033 &0.0037 &0.0033 \\
[3pt]
&A.SD &$\kappa$& &0.1181 &0.1105 & &0.0592 &0.0697 \\
&&$\alpha$& &0.0047 &0.0048 & &0.0027 &0.0028 \\
& &$\sigma$& &0.0030 &0.0030 & &0.0034 &0.0033 \\
[3pt]
&RMSD &$\kappa$& &0.0173 &0.0068 & &0.0424 &0.0100 \\
&&$\alpha$& &0.0004 &0.0001 & &0.0020 &0.0001 \\
& &$\sigma$& &0.0005 &0.0003 & &0.0027 &0.0001 \\
\hline
\end{tabular*}
\end{table}

The second part of the simulation was devoted to the diminishing
$\delta$ case. Here we wanted to confirm the differential behavior
of the AMLEs in the limiting distribution between $J=1$ and $J
\geq2$, as revealed in Section \ref{sec5}. The Vasicek model with
$\theta=(\kappa,\alpha,\sigma)'=(0.892,0.09,0.1817)'$ was
considered. We tried to create two scenarios: (i)
$n\delta^3\rightarrow\infty$ and (ii) $n\delta^3\rightarrow0$,
while $\delta\rightarrow0$. They were created by choosing
$\delta=n^{-1/6}$ and $\delta=n^{-1/2}$, respectively, whiling
selecting $n=500, 1\mbox{,}000, 2\mbox{,}000, 4\mbox{,}000$ and 8,000, respectively, to
create two streams of asymptotic sequences. For each $n$ and
$\delta$, we generated repeatedly the Vasicek sample paths 1,000
times. For each simulated sample path, we obtained the AMLEs
$\hat{\theta}_{n}^{(J)}$ for $J=1$ and $2$, respectively, and computed
the Wald statistics
\[
W_n(J)=n\bigl({\hat{\theta}_{n}^{(J)}}-\theta_0\bigr)'I(\delta)\bigl({\hat{\theta
}_{n}^{(J)}}-\theta_0\bigr).
\]

If $\sqrt{n}I^{1/2}(\delta)({\hat{\theta}_{n}^{(J)}}-\theta_0)$ is
asymptotically standard normally distributed in~$\mathbb{R}^d$,
then the Wald statistic $W_n(J) \stackrel{d}{\rightarrow} \chi_3^2$.
%$n(\hat{\theta}_n^{(1)}-\theta_0)^TI(\delta)(\hat{\theta}_n^{(1)}-
%will be asymptotic $\chi^2(3)$ distributed.
Based on the 1,000 Wald statistics from the simulations, we then
performed the Kolmogorov--Smirnov \mbox{(K--S)} test to test $H_0\dvtx W_n(J)
\sim\chi_3^2$, or not,
for each of the designed sequences of~$(n, \delta)$\vadjust{\goodbreak} generated under
the two scenarios.
%For the second situation, we picked $\delta=n^{-1/2}$ and the sample
%size were the same as the first situation. In this situation, we also
%used the %same way mentioned in the first situation to get the Wald
%statistic and performed Kolmogorov--Smirnov test to check whether $W_n
Table \ref{table4} reports the $p$-values of the test, which show that for $J=1$,
under both scenarios, the $p$-values of the K--S test became smaller, and
hence the above null hypothesis was rejected as $n$ increased. For
$J=2$, the $p$-values of the K--S test were sharply {different} between
the two scenarios.
%
%t4 ###
%t4 #&#
\begin{table}[b]
\caption{$p$-values of Kolmogorov--Smirnov test for
$W_n(J)\sim\chi^2_3$}
\label{table4}
\begin{tabular*}{\tablewidth}{@{\extracolsep{\fill}}lccd{1.5}d{1.5}@{}}
\hline
% after \ \hline or \cline{col1-col2} \cline{col3-col4} ...
\textbf{Situation} &$\bolds{n}$ & $\bolds{\delta}$
& \multicolumn{1}{c}{$\bolds{J=1}$}
& \multicolumn{1}{c@{}}{$\bolds{J=2}$} \\
\hline
$\delta=n^{-1/6}$& \hphantom{0,}500 & 0.3550 & 0.3524 & 0.0587 \\
&1,000 & 0.3162 & 0.4595 & 0.5830 \\
&2,000 & 0.2817 & 0.1149 & 0.2710 \\
&4,000 & 0.2510 & 0.0019 & 0.8309\\
&8,000 & 0.2236 & 5.74\mbox{e--}8 & 0.6002 \\
[4pt]
$\delta=n^{-1/2}$&\hphantom{0,}500 & 0.0447 & 5.04\mbox{e--}7 & 2.45\mbox{e--}8 \\
&1,000 & 0.0316 & 0.0003 & 9.72\mbox{e--}5\\
&2,000 & 0.0224 & 0.0006 & 0.0003\\
&4,000 & 0.0158 & 0.1109 & 0.0851\\
&8,000 & 0.0112 & 0.0470 & 0.0367\\
\hline
\end{tabular*}
\end{table}
In particular, the $p$-values were mostly quite large under the
scenario of $n \delta^3 \to\infty$, and they were\vspace*{1pt} largely significant
(small) when $\delta$ was diminishing at the faster rate of $n^{-1/2}$
such that $n\delta^3\rightarrow0$. These were consistent with our
theoretical findings in Section \ref{sec5}.

%apA #&#
\begin{appendix}\label{app}
\section*{Appendix}

We need the following technical assumptions in our analysis.

(A.1) (i) $\Theta$ is a compact set in $\mathbb{R}^d$, and the true
parameter $\theta_0$ is an interior point of $\Theta$; (ii) for all
values of the parameters $\theta$, Assumption 1--3 in
\citet{r2}  hold; (iii) the drift function $\mu(x;
\theta)$ is a bona fide function of~$\theta$ for each $x$.

(A.2) (i) For every $\delta> 0$, $ \mathbb{E}\nabla_\theta\log
f(X_t|X_{t-1},\delta;\theta_0)=0,%\mbox{and}\mathbb{E}\{\frac{\partial^2
%f(X_t|X_{t-1},\delta;\theta_0)}{\partial\theta_i\partial
$ and $\theta_0$ is the only root of $\mathbb{E}\nabla_\theta\log
f(X_t|X_{t-1},\delta;\theta)=0$. (ii) the MLE $\hat{\theta}_{n}$
and the $J$-term approximate MLE $ \hat{\theta}^{(J)}_{n}$ satisfy,
respectively,
\[
\sum_{t=1}^n\nabla_\theta\log
f(X_t|X_{t-1},\delta;{\hat{\theta}_{n}})=0
\quad\mbox{and}\quad\sum_{t=1}^n\nabla_\theta\log
f^{(J)}\bigl(X_t|X_{t-1},\delta;{\hat{\theta}^ {(J)}_{n}}\bigr)=0.
\]

%And (iii) $\bf\hat{\theta}_{n}$ is consistent to $\theta_0$
%and asymptotically normal such that (\ref{eqasynor}) is satisfied.

(A.3) There exist finite positive constants $\Delta$ and $K_1$
such that, for $l=1,2,3$, any $\delta\in(0,\Delta]$, $i_1, i_2, i_3
\in\{1,\ldots,d\}$ and $j=1$ and $2$,
\[
\mathbb{E}\sup_{\theta\in\Theta}\biggl\{\biggl|\frac{\partial^l
A_j(X_t|X_{t-1},\delta;\theta)}{\partial\theta_{i_1} \cdots
\partial\theta_{i_l} }\biggr|^3\biggr\}\leq
K_1.
\]

(A.4) There exist finite positive constants $\nu_l$ for $l=0, 1,2$
and $3$, $\Delta>0$ and $K_2$ such that $\nu_0 > 3$,
$\nu_2>\nu_1>3$, $\nu_3 > 1$ and
for any $i_1, \ldots, i_3 \in\{1,\ldots,d\}$ and $\delta\in(0,\Delta]$,
\[
\mathbb{E}\Biggl\{\sup_{\theta\in\Theta}\Biggl[\sum_{l=0}^\infty\biggl|\frac{\partial^q
c_l(\gamma(X_t;\theta)|\gamma(X_{t-1};\theta);\theta)}{\partial\theta_{i_1}
\cdots\partial\theta_{i_q}
}\biggr|\frac{\Delta^l}{l!}\Biggr]^{\nu_l}\Biggr\}\leq K_2.
\]

(A.5) For any $\delta>0$, the Fisher information matrix
\[
I(\delta):=-\mathbb{E}\nabla^2_{\theta\theta}\log
f(X_t|X_{t-1},\delta;\theta_0)
\]
is invertible and as $\delta\rightarrow0$ the largest eigenvalues
of $\delta I^{-1}(\delta)$ is bounded away from infinity.

(A.6) For each positive integer $K$, which may be infinite, and
any $\delta\in(0,\Delta]$,
\[
\mathbb{P}\Biggl\{\inf_{\theta\in\Theta}\Biggl|\sum_{l=0}^K
c_l(\gamma(X_t;\theta)|\gamma(X_{t-1};\theta);\theta)\frac{\delta
^l}{l!}\Biggr|=0\Biggr\}=0.
\]

(A.7) For any $\beta>1$ and $\eta>0$, there exists
$\Delta(\beta,\eta)>0$, then for any
$\delta\in(0,\Delta(\beta,\eta)]$ and $K$, where $K$ may be
infinite,
\[
\mathbb{P}\Biggl\{\inf_{\theta\in\Theta}\Biggl|\sum_{l=0}^Kc_l(\gamma(X_t;\theta
)|\gamma(X_{t-1};\theta);\theta)\frac{\delta^l}{l!}\Biggr|<\eta^{1/\beta}\Biggr\}
<\eta.
\]

Assumptions (A.1) and (A.2) are standard requirements for maximum
likelihood estimators. In particular, (A.1) (ii) contains conditions on
the smoothness of the drift and the diffusion which ensures the
existence of a~unique solution to (\ref{eq11}) as well as the infinite
differentiability of the transition density $f(x|x_0, \delta; \theta)$
with respect to $x$, $x_0$ and $\delta$, and three {times}
differentiable with respect to $\theta$ [\citet{r22}]. The second part of
(A.2) is the simplified approach of \citet{r16} assuming the MLEs
are the solutions of the likelihood score equations.
%We actually only need $f(X_{t+\delta}
%|X_t,\delta; \theta)$ to be continuously second differentiable with
%respect to $\theta$ almost surely. (A.2) is a regular assumption for
%the asymptotic normality in traditional MLE.
Assumption (A.3) is needed to guarantee the third derivative of $\log
f(X_t|X_{t-1},\delta;\theta)$ with respect to $\theta$ can be
controlled by an integrable function, while condition (A.4)
ensures the absolute convergence of the infinite series
${\sum_{l=0}^\infty}|c_l(\gamma(X_t;\theta)|\gamma(X_{t-1};\theta)|\delta^l/l!
= \exp\{\tilde{A}_3(x|x_0, \delta;\theta)\}$ as
\citet{r2}  has provided conditions on the nondegeneracy
of the
diffusion function and the boundary condition, which together with
the third part of condition (A.1) leads to the convergence of the
above infinite series $\exp\{\tilde{A}_3(x|x_0,
\delta;\theta)\}$. Condition (A.4) is also needed to allow exchange of
differentiation and summation for the infinite series.
The first part of the (A.5) is of standard in likelihood inference.
Its second part reflects the fact that for some processes
$\lim_{\delta\to0}I(\delta)$ may be singular, as conveyed in our
discussion in Section \ref{sec6} for the Vasicek process.
Condition~(A.6) is needed to guarantee the derivatives of log
transition density and log approximate transition density
exist with probability one. Condition~(A.7) is needed to manage the
denominators in the derivatives of the log approximate transition density,
ensuring that the probability of their taking small values can be
controlled uniformly.
% We note in particular, for any $K$, if $\delta=0$, then %$\inf_{
%We use this condition to control uniformly the probability for this
%class of random variables taking small values. }

We shall give the proofs for the propositions and theorems mentioned
in Sections \ref{sec3}--\ref{sec6}. We first present some lemmas about the true
transition density and its approximations, which we will use in
later proofs. The proofs for the lemmas can be found in \citet{r14}.
%
%le1 #&#
\begin{la}\label{lem1}
Under \textup{(A.1)} and \textup{(A.4)}, for any $\delta\in(0,\Delta)$, the infinite
series
\[
\sum_{l=0}^\infty
c_l(\gamma(X_t;\theta)|\gamma(X_{t-1};\theta))\frac{\delta^l}{l!}
\]
absolutely converges with probability 1, and for $k=1,2$ and $3$,
and ${i_1}, i_2, i_3 \in\{1,\ldots,d\}$,
\begin{eqnarray*}
&& \frac{\partial^k}{\partial\theta_{i_1} \cdots\partial
\theta_{i_k}}\sum_{l=0}^\infty
c_l(\gamma(X_t;\theta)|\gamma(X_{t-1};\theta))\frac{\delta^l}{l!} \\
&&\qquad= \sum_{l=0}^\infty
\frac{\partial^k}{\partial\theta_{i_1} \cdots\partial\theta_{i_k}}
c_l(\gamma(X_t;\theta)|\gamma(X_{t-1};\theta))\frac{\delta^l}{l!}.
\end{eqnarray*}
%
% \[
%c_l(\gamma(X_t;\theta)|\gamma(X_{t-1};\theta))\frac{\delta^l}{l!}=\sum
%_{l=0}^\infty\frac{\partial^2}%{\partial\theta_i\partial\theta_j}c_l(
% \]
% \[
%c_l(\gamma(X_t;\theta)|\gamma(X_{t-1};\theta))\frac{\delta^l}{l!}=\sum
%_{l=0}^\infty\frac{\partial^3}%{\partial\theta_i\partial\theta_j
% \]
\end{la}

%convergence of the infinite series.
%%\[
%%\sum_{l=0}^\infty
%%c_l(\gamma(X_t;\theta)|\gamma(X_{t-1};\theta))\frac{\delta^l}{l!}.
%%\]
%Let
%$S_n(\delta)=\sum_{l=0}^nc_l(\gamma(X_t;\theta)|\gamma(X_{t-1};\theta))
%
%For a fixed $\delta\in(0,\Delta)$ and $\theta\in\Theta$,
%N}|S_m(\delta)-S_M(\delta)|>\varepsilon\}\leq\mathbb{P}\{\sum
%_{l=M+1}^N|c_l(\gamma(X_t;\theta)|\gamma(X_{t-1};\theta))|\frac{
%Applying Markov inequality,
%N}|S_m(\delta)-S_M(\delta)|>\varepsilon\}\leq\varepsilon^{-2}\cdot
%Letting $N\rightarrow\infty$, we get from (A.4),
%M}|S_m-S_M|>\varepsilon\}\leq\varepsilon^{-2}\cdot\mathbb{E}\{\sum
%_{l=M+1}^\infty|c_l(\gamma(X_t;\theta)|\gamma(X_{t-1};\theta))|\frac{
%
%If we let $\omega_M=\sup_{m,n\geq M}|S_m-S_n|$, then
%$\omega_M\downarrow$ as $M\uparrow$ and
%M}|S_m-S_M|>\varepsilon\}\rightarrow0
%as $M\rightarrow\infty$. Hence, $\omega_M\downarrow0$ almost surely.
%Then, we attain the absolutely convergence of the infinite series.
%Actually, this absolute convergence is uniform on $\Theta$.
%
%Next, we consider the exchange between the differentiation and the
%summation. The key is to prove that
%is uniformly convergent on $\Theta$ with probability 1. Using the
%same method above and from (A.4), the result is correct. Then,
%c_l(\gamma(X_t;\theta)|\gamma(X_{t-1};\theta))\frac{\delta^l}{l!}=\sum
%_{l=0}^\infty\frac{\partial}{\partial\theta_i}c_l(\gamma(X_t;\theta)|
%for any $i\in\{1,\cdots,d\}$ with probability 1. Using the same
%approach, we can show the exchange between differentiation and the
%summation is also valid for $k=2$ and $3$, respectively.
%
%
%le2 #&#
\begin{la}\label{lem2}
Under \textup{(A.6)} and \textup{(A.7)}, for any positive $\beta>1$, there exist two
constants $m(\beta)<\infty$ and $\Delta_1(\beta)>0$ such that for
any $\delta\in(0,\Delta_1(\beta)]$ and~$J$, where~$J$ can be
infinity, then
\[
\mathbb{E}\Biggl\{\sup_{\theta\in\Theta}\Biggl|\sum_{l=0}^J
c_l(\gamma(X_t;\theta)|\gamma(X_{t-1};\theta))
\frac{\delta^l}{l!}\Biggr|^{-\beta}\Biggr\}<m(\beta).
\]
\end{la}

\begin{la}\label{lem3}
Under \textup{(A.1), (A.3), (A.4), (A.6), (A.7)}, there exist two constants
$M_1<\infty$ and $\Delta_2>0$ such that, for any $J$, where $J$ can
be infinity, $\delta\in(0,\Delta_2)$ and $i,j,k\in\{1,\ldots,d\}$,
\[
\mathbb{E}\biggl\{\sup_{\theta\in\Theta}\biggl|\frac{\partial^3\log
f^{(J)}(X_t|X_{t-1},\delta;\theta)}{\partial\theta_i\,\partial\theta
_j\,\partial\theta_k}\biggr|\biggr\}<M_1.
\]
\end{la}

\begin{pf*}{Proof of Proposition \ref{prop1}}
Using the same method in the proof
of Lemma~\ref{lem3}, we know (a) holds. On the other hand, Lemma \ref{lem3}
implies
(b).
\end{pf*}
\begin{pf*}{Proof of Proposition \ref{prop2}} %By the definition of
See the proof of Proposition 2 in \citet{r14}.
\end{pf*}
\begin{pf*}{Proof of Proposition \ref{prop3}}
Recall Proposition \ref{prop2}, then
\[
\|I^{-1}(\delta) N(\theta_0, J,
\delta)+E_d\|_2\leq\|I^{-1}(\delta)\|_2\cdot\| N(\theta_0, J,
\delta)+I(\delta)\|_2\leq C\delta^{J}.
\]
If $C\delta^{J}<1$, then
\[
\|N^{-1}(\theta_0,J,\delta)I(\delta)+E_d\|_2\leq\frac
{\|I^{-1}(\delta)
N(\theta_0, J, \delta)+E_d\|_2}{1-\|I^{-1}(\delta) N(\theta_0, J,
\delta)+E_d\|_2}.
\]

From Proposition \ref{prop2}, if $C\delta^{J+1}<1$, then
\[
\|N^{-1}(\theta_0,J,\delta)+I^{-1}(\delta)\|_2\leq\frac
{\|I^{-1}(\delta)\|_2^2\|N(\theta_0,J,\delta)+I(\delta
)\|_2}{1-\|I^{-1}(\delta)\|_2\|N(\theta_0,J,\delta)+I(\delta)\|_2}.
\]

On the other hand, using the same method in the proof of Proposition
\ref{prop2}, we have
\[
\|U(\theta_0,J,\delta)\|_2\leq C\delta^{J+1}
\]
for any positive $J$ and $\delta\in(0,\bar{\Delta})$. Hence, we can
find the constants $C_1, C_2$ and $\underline{\Delta}>0$ such that
\[
\|N^{-1}(\theta_0,J,\delta)I(\delta)+E_d\|_2\leq
C_1\delta^{J} \quad\mbox{and}\quad \|N^{-1}(\theta_0,J,\delta)U(\theta_0,J,\delta
)\|_2\leq
C_2\delta^{J}
\]
for any positive integer $J$ and $\delta\in(0,\underline{\Delta})$.
\end{pf*}
\begin{pf*}{Proof of Proposition \ref{prop4}} Use the same method in
the proof of Proposition~\ref{prop2}.
\end{pf*}
\begin{pf*}{Proof of Proposition \ref{prop5}} We'll use Corollary 2.1
in Newey (\citeyear{r28}) to prove this proposition. We only need to verify three
conditions under two situations mentioned in Proposition \ref{prop5}:

\begin{longlist}
\item
for any $i\in\{1,\ldots,d\}$,
\[
\mathbb{E}\biggl\{\frac{\partial}{\partial\theta_i}\log
f(X_t|X_{t-1},\delta;\theta)\biggr\}\qquad \mbox{is equicontinuous};
\]

\item for any $i\in\{1,\ldots,d\}$,
\[
\sup_{\theta\in\Theta}\Biggl\|\frac{1}{n}\sum_{t=1}^n\frac{\partial
^2}{\partial\theta_i\,\partial\theta'}\log
f(X_t|X_{t-1},\delta;\theta)\Biggr\|_2=O_p(1);
\]

\item for any $i\in\{1,\ldots,d\}$ and $\theta\in\Theta$,
\[
\frac{1}{n}\sum_{t=1}^n\frac{\partial}{\partial\theta_i}\log
f(X_t|X_{t-1},\delta;\theta)-\mathbb{E}\biggl\{\frac{\partial}{\partial\theta
_i}\log
f(X_t|X_{t-1},\delta;\theta)\biggr\}\stackrel{p}{\rightarrow}0.
\]
\end{longlist}
For any $\theta^*,\theta^{**}\in\Theta$, note that
\begin{eqnarray*}
&& \mathbb{E}\biggl\{\frac{\partial}{\partial\theta_i}\log
f(X_t|X_{t-1},\delta;\theta^*)\biggr\}-\mathbb{E}\biggl\{\frac{\partial}{\partial
\theta_i}\log
f(X_t|X_{t-1},\delta;\theta^{**})\biggr\}\\
&&\qquad =\mathbb{E}\biggl\{\frac{\partial^2}{\partial\theta_i\,\partial\theta'}\log
f(X_t|X_{t-1},\delta;\bar{\theta})\biggr\}\cdot(\theta^*-\theta^{**}),
\end{eqnarray*}
where $\bar{\theta}$ is on the joint line between $\theta^*$ and
$\theta^{**}$. Then
\begin{eqnarray*}
&& \biggl|\mathbb{E}\biggl\{\frac{\partial}{\partial\theta_i}\log
f(X_t|X_{t-1},\delta;\theta^*)\biggr\}-\mathbb{E}\biggl\{\frac{\partial}{\partial
\theta_i}\log
f(X_t|X_{t-1},\delta;\theta^{**})\biggr\}\biggr|\\
&&\qquad \leq\biggl\|\mathbb{E}\biggl\{\frac{\partial^2}{\partial\theta_i\,\partial
\theta'}\log
f(X_t|X_{t-1},\delta;\bar{\theta})\biggr\}\biggr\|_2
\cdot\|\theta^{*}-\theta^{**}\|_2.
\end{eqnarray*}
For any $j\in\{1,\ldots,d\}$, use the same method in the proof of
Lemma \ref{lem3}, we know that there exists a constant $C$, which is not
dependent on $J$ and $\delta$, and $\hat{\Delta}>0$ such that, for
any $J$ and $\delta\in(0,\hat{\Delta}]$,
\[
\mathbb{E}\biggl\{\sup_{\theta\in\Theta}\biggl|\frac{\partial^2}
{\partial\theta
_i\,\partial\theta_j}\log
f(X_t|X_{t-1},\delta;\theta)\biggr|\biggr\}<C.
\]
Hence, (i) and (ii) can be established.

To verify (iii), from (A.3) [Lemmas 3 and 4 in
\citet{r7}], we know that there exists a positive constant
$\kappa$ such that for any $t_1<t_2$,
\begin{eqnarray*}
& &\biggl|\mathbb{E}\biggl\{\biggl[\frac{\partial}{\partial\theta_i}\log
f(X_{t_1}|X_{t_1-1},\delta;\theta)-\mathbb{E}\biggl\{\frac{\partial}{\partial
\theta_i}\log
f(X_{t_1}|X_{t_1-1},\delta;\theta)\biggr\}\biggr]\\
&&\hspace*{6pt}\quad \times\biggl[\frac{\partial}{\partial\theta_i}\log
f(X_{t_2}|X_{t_2-1},\delta;\theta)-\mathbb{E}\biggl\{\frac{\partial}{\partial
\theta_i}\log
f(X_{t_2}|X_{t_2-1},\delta;\theta)\biggr\}\biggr]\biggr\}\biggr|\\
&&\hspace*{6pt}\qquad \leq C\cdot\exp\{-\kappa(t_2-t_1)\delta\},
\end{eqnarray*}
where
\[
C=\mathbb{E}\biggl\{\biggl[\frac{\partial}{\partial\theta_i}\log
f(X_{t}|X_{t-1},\delta;\theta)-\mathbb{E}\biggl\{\frac{\partial}{\partial
\theta_i}\log
f(X_{t}|X_{t-1},\delta;\theta)\biggr\}\biggr]^2\biggr\}.
\]
Then
\begin{eqnarray*}
&& \mathbb{E}\Biggl\{\frac{1}{n}\sum_{t=1}^n\biggl[\frac{\partial}{\partial\theta
_i}\log
f(X_{t}|X_{t-1},\delta;\theta)-\mathbb{E}\biggl\{\frac{\partial}{\partial
\theta_i}\log
f(X_{t}|X_{t-1},\delta;\theta)\biggr\}\biggr]\Biggr\}^2\\
&&\qquad \leq\frac{C}{n}+\frac{C}{n}\cdot\frac{\exp\{-\kappa\delta\}
}{1-\exp\{-\kappa\delta\}}\\
&&\qquad \leq3\biggl[2K_1+K_2\cdot
m\biggl(\frac{2\nu_1}{\nu_1-2}\biggr)\biggr]
\cdot\biggl\{\frac{1}{n}+\frac{1}{n[\exp(\kappa
\delta)-1]}\biggr\}\rightarrow0,
\end{eqnarray*}
under the two situations mentioned in the statement of Proposition
\ref{prop5}. Hence we complete the proof.
\end{pf*}
\begin{pf*}{Proof of Proposition \ref{prop6}} From (A.2), we can get
$n^{-1}\nabla_\theta\ell_{n,\delta}(\hat{\theta}_n)=0$. Expanding it
at $\theta_0$,
\[
0=\frac{1}{n}\sum_{t=1}^n\nabla_\theta\log
f(X_t|X_{t-1},\delta;\theta_0)+\frac{1}{n}\sum_{t=1}^n\nabla^2_{\theta
\theta}\log
f(X_t|X_{t-1},\delta;\tilde{\theta})\cdot(\hat{\theta}_n-\theta_0).
\]
Then
\[
\hat{\theta}_n-\theta_0=\Biggl\{-\frac{1}{n}\sum_{t=1}^n\nabla^2_{\theta\theta
}\log
f(X_t|X_{t-1},\delta;\tilde{\theta})\Biggr\}^{-1}\cdot\frac{1}{n}\sum
_{t=1}^n\nabla_\theta
\log f(X_t|X_{t-1},\delta;\theta_0).
\]
Define\vspace*{1pt} $I_n(\delta)=-n^{-1}\sum_{t=1}^n\nabla^2_{\theta\theta}\log
f(X_t|X_{t-1},\delta;\theta_0)$. From Lemma \ref{lem3}, an\break
$-{n}^{-1}\times\sum_{t=1}^n\nabla^2_{\theta\theta}\log
f(X_t|X_{t-1},\delta;\tilde{\theta})=I_n(\delta)\cdot\{1+o_p(1)\}$.
Using the same way as that in the verification of (iii) in the proof
of Proposition \ref{prop5}, we can get
$I_n(\delta)-I(\delta)=O_p\{(n\delta)^{-1/2}\}$. If
$n\delta^{3}\rightarrow\infty$, by (A.5),
\begin{eqnarray*}
\Biggl\{-\frac{1}{n}\sum_{t=1}^n\nabla^2_{\theta\theta}\log
f(X_t|X_{t-1},\delta;\tilde{\theta})\Biggr\}^{-1}
&=&\bigl\{I(\delta)\cdot\{1+o_p(1)\}+O_p\{(n\delta)^{-1/2}\}\bigr\}
^{-1}\\
&=&I^{-1}(\delta)\cdot\{1+o_p(1)\}.
\end{eqnarray*}
%
%Then,
Then
\[
\sqrt{n}I^{1/2}(\delta)(\hat{\theta}_n-\theta_0)=I^{-1/2}(\delta)\frac
{1}{n^{1/2}}\sum_{t=1}^n\nabla_\theta
\log f(X_t|X_{t-1},\delta;\theta_0)\cdot\{1+o_p(1)\}.
\]
We will use the martingale central limit theorem [\citet{r12},
page~476] to show that the first part on the right-hand side of the
above equation
converges to a standard normal distribution. For any\vspace*{1pt}
$\alpha\in\mathbb{R}^d$ with unit $L_2$ norm, to simplify notations,
let $ U_{n,m}=\alpha'I^{-1/2}(\delta)n^{-1/2}\nabla_\theta\log
f(X_m|X_{m-1},\delta;\theta_0)$ and
$\mathscr{F}_{n,m}=\sigma(X_1,\ldots,X_m)$. It is easy to check
$(U_{n,m},\mathscr{F}_{n,m})$ is a martinga\-le~difference array. By the
Markov property and Birkhoff's Ergodic theo\-rem,~%
$V_{n,n}=\sum_{m=1}^n\mathbb{E}(U_{n,m}^2|\mathscr{F}_{n,m})
\stackrel{p}{\rightarrow}\mathbb{E}U_{n,m}^2=1$.
On the other hand,\break $\sum_{m=1}^n|U_{n,m}|^3\leq
C(n\times\delta^3)^{-1/2}\rightarrow0$. This implies the asymptotic
normality of
$\sqrt{n}\alpha'I^{1/2}(\delta)(\hat{\theta}_n-\theta_0)$. Then we
complete the proof.
%Note the
%arbitrary selection of $\alpha$, we complete the proof.
\end{pf*}
\begin{pf*}{Proof of Theorem \ref{theo1}}
From Propositions \ref{prop4} and \ref{prop5}, we can get
\[
\bigl\|\mathbb{E}\nabla_\theta\log
f\bigl(X_t|X_{t-1},\delta;{\hat{\theta}_{n}^{(J)}}\bigr)\bigr\|_2
%f(X_t|X_{t-1},\delta;\theta)-\mathbb{E}\{\frac{\partial}{\partial
%f(X_t|X_{t-1},\delta;\theta)\}{|}{|}_2\\
%&+\sup_{\theta\in\Theta}{|}{|}\frac{1}{n}\sum_{t=1}^n\frac{\partial}{
%f(X_t|X_{t-1},\delta;\theta)-\frac{1}{n}\sum_{t=1}^n\frac{\partial}{
%f^{(J)}(X_t|X_{t-1},\delta;\theta){|}{|}_2\\
\stackrel{p}{\rightarrow}0
\]
for either: (i) $\delta\in(0,\tilde{\Delta}\wedge\dot{\Delta}]$
being fixed, $J\rightarrow\infty$ and $n\rightarrow\infty$, or (ii)
$J$ being fixed, $n\rightarrow\infty$, $\delta\rightarrow0$ but
$n\delta\rightarrow\infty$. Hence, noting condition (A.2)(i), we
have the consistency of the AMLE $\hat{\theta}_{n}^{(J)}$.
\end{pf*}
\begin{pf*}{Proof of Theorem \ref{theo2}}
For fixed $\delta$, from Theorem \ref{theo1} and
(\ref{eq41}), we know that the leading order term of
$\hat{\theta}_{n}^{(J)}-\theta_0$ contains two parts: one is
$N^{-1}U_n$, and the other is
$N^{-1}(N_n+F_n)({\hat{\theta}_{n}}-\theta_0)$. Hence,
${\hat{\theta}_{n}^{(J)}}-\theta_0=O_p\{\delta^{J+1}+(n\delta)^{-1/2}\}$.

For $J$ fixed and $\delta\rightarrow0$, Proposition \ref{prop4} implies
%f(X_t|X_{t-1},\delta;{\bf\hat{\theta}_{n}^{(J)}}){|}{|}_2\}\leq
%C\delta^{J+1}.
%Then,
%
\begin{eqnarray*}
&& \mathbb{E}\Biggl\{\Biggl\|\frac{1}{n}\sum_{t=1}^n\nabla_\theta\log
f\bigl(X_t|X_{t-1},\delta;{\hat{\theta}_{n}^{(J)}}\bigr)-\frac{1}{n}\sum
_{t=1}^n\nabla_\theta\log
f(X_t|X_{t-1},\delta;{\hat{\theta}_{n}})\Biggr\|_2\Biggr\}\\
&&\qquad \leq C\delta^{J+1}.
\end{eqnarray*}
This means that
\[
\mathbb{E}\Biggl\{\Biggl\|\frac{1}{n}\sum_{t=1}^n\nabla^2_{\theta\theta}\log
f(X_t|X_{t-1},\delta;\tilde{\theta})\cdot\bigl({\hat{\theta}_{n}^{(J)}-\hat
{\theta}_{n}}\bigr)\Biggr\|_2\Biggr\}\leq
C\delta^{J+1},
\]
where $\tilde{\theta}$ is on the joining line between
$\hat{\theta}_{n}^{(J)}$ and $\hat{\theta}_{n}$. Hence
\[
\frac{1}{n}\sum_{t=1}^n\nabla^2_{\theta\theta}\log
f(X_t|X_{t-1},\delta;\tilde{\theta})\cdot\bigl({\hat{\theta}_{n}^{(J)}-\hat
{\theta}_{n}}\bigr)=O_p(\delta^{J+1}).
\]
Since $\tilde{\theta}\stackrel{p}{\rightarrow}\theta_0$ and
${\hat{\theta}_{n}^{(J)}-\hat{\theta}_{n}}=o_p(1)$,
\[
\frac{1}{n}\sum_{t=1}^n\nabla^2_{\theta\theta}\log
f(X_t|X_{t-1},\delta;\theta_0)\cdot\bigl({\hat{\theta}_{n}^{(J)}-\hat{\theta
}_{n}}\bigr)=O_p(\delta^{J+1}).
\]
On the other hand, from Proposition \ref{prop2}, we know
\begin{eqnarray*}
&&\frac{1}{n}\sum_{t=1}^n\nabla^2_{\theta\theta}\log
f(X_t|X_{t-1},\delta;\theta_0)-\frac{1}{n}\sum_{t=1}^n\nabla^2_{\theta
\theta}\log
f^{(J)}(X_t|X_{t-1},\delta;\theta_0)\\
&&\qquad=O_p(\delta^{J+1}).
\end{eqnarray*}
Then\vspace*{1pt}
$N_n({\hat{\theta}_{n}^{(J)}-\hat{\theta}_{n}})=O_p(\delta^{J+1})$.
Using the same way of verifying (iii) in the proof of Proposition \ref{prop5},
we know $N_n-N=O_p\{(n\delta)^{-1/2}\}$. As
$n\delta^3\rightarrow\infty$, then
$N({\hat{\theta}_{n}^{(J)}-\hat{\theta}_{n}})=O_p(\delta^{J+1})$.
Hence, $ {\hat{\theta}_{n}^{(J)}-\hat{\theta}_{n}}=O_p(\delta^J). $
At the same time, we know
${\hat{\theta}_{n}}-\theta_0=O_p\{(n\delta)^{-1/2}\}$. Then
\[
{\hat{\theta}_{n}^{(J)}}-\theta_0=O_p\{\delta^J+(n\delta)^{-1/2}\}.
\]
This completes the proof of Theorem \ref{theo2}.
\end{pf*}
\begin{pf*}{Proof of Theorem \ref{theo4}} We only need to prove following result:
\[
\sqrt{n}\tilde{I}_n^{1/2}\bigl(\hat{\theta}_n^{(J)},J,\delta\bigr)\bigl(\hat{\theta
}_n^{(J)}-\theta_0\bigr)=\sqrt{n}I^{1/2}(\delta)\bigl(\hat{\theta}_n^{(J)}-\theta
_0\bigr)+o_p(1)
\]
under\vspace*{1pt} the two situations mentioned in Theorem \ref{theo4}. Using the approach
in the proof of Lemma \ref{lem3}, we have
$\tilde{I}_n(\hat{\theta}_n^{(J)},J,\delta)-\tilde{I}_n({\theta
}_0,J,\delta)=O_p\{\|\hat{\theta}_n^{(J)}-\theta_0\|_2\}$.
Also,\vspace*{1pt} using the same way of verifying (iii)\vspace*{1pt} in the proof of
Proposition \ref{prop5},
$\tilde{I}_n({\theta}_0,J,\delta)-\mathbb{E}\tilde{I}_n({\theta
}_0,J,\delta)=O_p\{(n\delta)^{-1/2}\}$.
By the same argument in the proof of Proposition~\ref{prop2},
$\mathbb{E}\tilde{I}_n({\theta}_0,J,\delta)-I(\delta)=O(\delta^{J+1})$.
Hence, if $n\delta^3\rightarrow\infty$, under either asymptotic
regime in Theorem \ref{theo4},
\[
\tilde{I}_n^{1/2}\bigl(\hat{\theta}_n^{(J)},J,\delta\bigr)=I^{1/2}(\delta
)\cdot\{1+o_p(1)\}.
\]
Then we complete the proof.%\vadjust{\goodbreak}
\end{pf*}
\end{appendix}

\section*{Acknowledgments}

We thank the Associate Editor for very constructive comments and
suggestions which have improved the presentation of the paper.

The first author thanks Department of Statistics at Iowa State
University for hospitality during his visits.

%suskaldyti doi

% imsref loaded by lrinkeviciute, 2011-11-23 15:51:36
% imsref loaded by lrinkeviciute, 2011-11-23 15:56:37

\printaddresses

\end{document}